\documentclass[11pt, a4paper]{amsart}
\usepackage{amsmath,amssymb,amsfonts, float,mathrsfs}
\usepackage{thmtools}

\usepackage[all]{xy}
\usepackage{caption}
\usepackage{enumerate}
\usepackage{a4wide}
\usepackage{diagbox}
\usepackage{scalerel}

\usepackage{color,graphicx}
\usepackage{pgfplots}
\usepackage{pgf}
\usepackage{tikz}
\usetikzlibrary{patterns}
\usepgfplotslibrary{patchplots} 
\usetikzlibrary{pgfplots.patchplots} 
\pgfplotsset{width=9cm,compat=1.5.1}

\usepackage{bm}
\usepackage{graphicx}
\usepackage[breaklinks=true]{hyperref}

\usepackage{xcolor}
\usepackage{multicol}
\usepackage{tabularx}

\usepackage{hyperref}
\usepackage[capitalize]{cleveref}

\usepackage[normalem]{ulem}

\newtheorem{thm}{Theorem}[section] 
\newtheorem*{thm*}{Theorem} 
\newtheorem{prop}[thm]{Proposition}
\newtheorem{lem}[thm]{Lemma}
\newtheorem{cor}[thm]{Corollary}

\theoremstyle{definition}
\newtheorem{definition}[thm]{Definition}
\newtheorem{expl}[thm]{Example}

\newtheorem{rem}[thm]{Remark}

\DeclareMathOperator{\R}{\mathbb{R}}
\DeclareMathOperator{\C}{\mathbb{C}}
\DeclareMathOperator{\Z}{\mathbb{Z}}
\DeclareMathOperator{\N}{\mathbb{N}}
\DeclareMathOperator{\F}{\mathbb{F}}

\DeclareMathOperator{\Q}{\mathbb{Q}}

\DeclareMathOperator{\Hom}{{\rm Hom}}

\DeclareMathOperator{\Res}{{\rm Res}}

\DeclareMathOperator{\OO}{\mathcal{O}}

\DeclareMathOperator{\GL}{\text{GL}}

\DeclareMathOperator{\one}{\mathbb{\mathbf{1}}}

\newcommand{\bc}{\mathbf c}
\newcommand{\be}{\mathbf e}

\newcommand{\x}{\mathbf x}

\newcommand{\bv}{\mathbf v}

\newcommand{\zero}{\mathbf 0}
\newcommand{\spec}{\operatorname{spec}}

\newcommand{\trivialS}{ {\rho_{0}}_{\scaleto{S}{5pt}}}

    \DeclareFontFamily{U}{wncy}{}
    \DeclareFontShape{U}{wncy}{m}{n}{<->wncyr10}{}
    \DeclareSymbolFont{mcy}{U}{wncy}{m}{n}
    \DeclareMathSymbol{\Sha}{\mathord}{mcy}{"58}

\numberwithin{equation}{section}

\newcommand{\Cay}{\textrm{Cay}}

\DeclareSymbolFont{bbold}{U}{bbold}{m}{n}
\DeclareSymbolFontAlphabet{\mathbbold}{bbold}

\usepackage{hyperref}

\begin{document}

\title{Sedentariness in Cayley graphs}

\author{Karen Meagher, Hermie Monterde, Tung T. Nguyen }
\address{Department of Mathematics and Statistics, University of Regina, SK, Canada S4S 0A2}
\email{karen.meagher@uregina.ca}

\address{Department of Mathematics and Statistics, University of Regina, SK, Canada S4S 0A2}
\email{hermie.monterde@uregina.ca}
 
\address{
Department of Mathematics, Elmhurst University} 
\email{tung.nguyen@elmhurst.edu}

\keywords{Cayley graph, Quantum walk, Sedentary vertex, Twins, Necklace polynomial}
\subjclass[2020]{Primary 05C25; Secondary 05C50, 81P45, 15A16, 15A18}

\vspace{-0.1in}
\begin{abstract}
Let $A$ be the adjacency matrix of a finite graph $\Gamma$. We say that vertex $u$ of $\Gamma$ is sedentary if for some constant $C>0$, $|\exp(itA)_{u,u}|\geq C$ for all $t\in\R$. In this paper, we investigate vertex sedentariness in walk-regular graphs with a focus on Cayley graphs and, more generally, vertex-transitive graphs. In particular, we define the necklace polynomial of an integral graph and use this notion to characterize sedentariness in integral walk-regular graphs. As an application, we use basic tools from algebraic number theory to show that an integral Cayley graph on a group of odd order is sedentary and provide a lower bound on its sedentariness. We also give sufficient conditions such that a unitary Cayley graph over a finite ring is sedentary. Furthermore, we characterize the existence of twin vertices in Cayley graphs and vertex-transitive graphs, and completely determine the occurrence of sedentariness, perfect state transfer and pretty good state transfer in vertex-transitive graphs with twins. Finally, we prove that if $\Gamma(G,S)$ is a nonempty, integral, normal Cayley graph, then for any prime $p$ with $\gcd(p, |G|)=1$, there is at least one eigenvalue of $\Gamma(G,S)$ which is not a multiple of $p$. 
This generalizes a result of Arnadottir and Godsil, which states that an integral normal Cayley graph of odd order possesses an odd eigenvalue.
\end{abstract}

\maketitle

\vspace{-0.3in}

\tableofcontents

\section{Introduction}

Throughout this paper, we let
$\Gamma$ be a graph on $n$ vertices with vertex set $V(\Gamma)$ and adjacency matrix $A=A(\Gamma)$. A \emph{(continuous) quantum walk} on $\Gamma$ is determined by a one-parameter family of complex symmetric $n\times n$ unitary matrices:
\begin{equation*}
    U(t)=e^{\rm{i}tA}\quad t\in\R.
\end{equation*}
We call $U(t)$ the \textit{transition matrix} of $\Gamma$. A quantum walk on $\Gamma$ describes the propagation of quantum states in the spin network modeled by $\Gamma$. Since $U(t)$ is unitary for all $t$, we have
\begin{equation}
\label{eq:unitary}
\sum_{j=1}^n|U(t)_{u,j}|^2=1.
\end{equation}
So, we may interpret $|U(t)_{u,v}|^2$ as the probability that the quantum state assigned at vertex $u$ is transmitted to vertex $v$ at time $t$. If $|U(\tau)_{u,v}|=1$ for some time $\tau>0$, then we say that \textit{perfect state transfer (PST)} occurs between $u$ and $v$ at time $\tau$. A relaxation of PST is called \textit{pretty good state transfer} (PGST), which occurs between $u$ and $v$ if there is a sequence $\{\tau_k\}$ such that $\lim_{k\rightarrow\infty}|U(\tau_k)_{u,v}|=1$. On the other hand, if there is a constant $C>0$ such that
\begin{equation}
\label{eq:sed}
\inf_{t>0}\big\lvert U(t)_{u,u}\big\rvert\geq C,
\end{equation}
then we say that vertex $u$ is \textit{$C$-sedentary}. If $C$ does not matter, then we simply say sedentary in  lieu of $C$-sedentary. From (\ref{eq:unitary}), it follows that PGST and vertex sedentariness are mutually exclusive. That is, a vertex involved in PGST is not sedentary, and a sedentary vertex is not involved in PGST. Further, we say that $\Gamma$ is \textit{sedentary} if each vertex of $\Gamma$ is sedentary.

For a sedentary vertex $u$, define $\mathscr{C}:=\inf_{t>0}\big\lvert U(t)_{u,u}\big\rvert$. Note that $\mathscr{C}$ is the largest possible $C>0$ such that vertex $u$ is $C$-sedentary. It is desirable to obtain the exact value, or a good approximation, of $\mathscr{C}$ whenever possible, for this allows us to provide a rough upper bound for $\big\lvert U(t)_{u,v}\big\rvert$ for all $t\in\R$. Indeed, for two vertices $u$ and $v$ in $\Gamma$, we have  $|U(t)_{u,u}|^2+|U(t)_{u,v}|^2\leq 1$ for all $t\in\R$ by (\ref{eq:unitary}). Thus, if $u$ is sedentary, then $|U(t)_{u,u}|\geq \mathscr{C}$ for all $t$, and so
\begin{equation*}
\big\lvert U(t)_{u,v}\big\rvert\leq \sqrt{1-\mathscr{C}^2}\quad \text{for all $t\in\R$}.
\end{equation*}
In \cite{cooper2026quantum}, it is shown that the above bound is related to peak state transfer, which is another relaxation of PST first studied in \cite{coutinho2025peak}.

This paper is organized as follows.  Section \ref{sec:pre} is allotted to basic definitions and notation relevant to the rest of the paper. In Section~\ref{sec:tw}, we recall known results about sedentariness in twins, and completely characterize PST, PGST and vertex sedentariness in graphs that are wreath products with a complete or empty graph (also known as blow-up graphs in the literature) (Theorem \ref{thm:sedwreath} and Corollary \ref{cor:blup}). In Section~\ref{sec:wr}, we provide sufficient conditions for sedentariness in walk-regular graphs (Corollary \ref{cor:sed}) and show that a vertex in a walk-regular graph is sedentary if and only if every vertex is (Proposition \ref{prop:sedwrg}). In Section~\ref{sec:integ}, we define the necklace polynomial of an integral graph and use this to characterize sedentariness in integral walk-regular graphs. In particular, we show that an integral walk-regular graph is sedentary if and only if its associated necklace polynomial has no roots on the unit circle (Theorem \ref{thm:rootsedentariness}). Using tools from basic algebraic number theory, we show that a connected, integral, walk-regular graph is sedentary whenever the multiplicity of each of its eigenvalues, except the valency, is a multiple of some integer greater than one (Theorem \ref{thm:main-multiple-eig}). Moreover, for this family of sedentary graphs, we use the theory of resultants to derive a lower bound on the sedentariness (Theorem \ref{thm:C}). In Section~\ref{sec:cayley}, we apply our results to Cayley graphs. In particular, we show that an integral Cayley graph on a group of odd order is sedentary (Theorem \ref{thm:oddsed}), and derive sufficient conditions such that a $U$-unitary Cayley graph is sedentary. Our goal in Section~\ref{sec:twcayley} is two-fold: we characterize the existence of twins in Cayley graphs and vertex-transitive graphs (Lemmas \ref{lem:conditions_for_twins}, \ref{lem:conditions_for_twins2}, and Theorem \ref{thm:conditions_for_twins_vert-trans}), and we use this characterization to completely determine PST, PGST and sedentariness in vertex-transitive graphs with twins (Theorem \ref{thm:caytwsed} and Corollary \ref{cor:pstvertextranstw}). In Section~\ref{sec:modprop}, we also investigate whether certain modular properties hold for eigenvalues of integral graphs and Cayley graphs. Finally, we present open questions in Section~\ref{sec:futwork}.

\medskip

\section{Preliminaries}\label{sec:pre}

Let $\Gamma$ be a graph on $n$ vertices. Throughout the paper, we refer to the eigenvalues and eigenvectors of $A(\Gamma)$ as the eigenvalues and eigenvectors of $\Gamma$. We order the eigenvalues of $\Gamma$ in a nonincreasing order $\lambda_1 \geq \lambda_2 \ge \cdots \geq \lambda_n$.  The \textit{spectrum} of $\Gamma$, denoted $\operatorname{spec}(\Gamma)$, is the multiset of eigenvalues of $A(\Gamma)$ (counted with multiplicity). For convenience, we also denote by $\sigma(G)$ the set of distinct eigenvalues of $\Gamma.$ We sometimes use $\spec(\Gamma)$ and $\sigma(\Gamma)$ interchangeably when multiplicities are not relevant. The meaning will be clear from the context. For each $\lambda\in\sigma(\Gamma)$, we let $E_\lambda$ denote the orthogonal projection matrix associated with the $\lambda$-eigenspace of $A(\Gamma)$. The \textit{(eigenvalue) support} of vertex $u$ in $\Gamma$, denoted $\Phi_{u}(\Gamma)$, is the set
\begin{equation*}
\Phi_{u}(\Gamma)=\{\lambda\in\sigma(\Gamma):(E_\lambda)_{u,u}\neq 0\}.
\end{equation*}
If $\Gamma$ is clear from the context, then we write $\Phi_{u}(\Gamma)$ as $\Phi_{u}$. We also remark that $\lambda\in\Phi_{u}(\Gamma)$ if and only if there is an eigenvector $\mathbf{w}$ for $\Gamma$ associated with $\lambda$ that has a nonzero $u$th entry.

We write $M^T$ and $M^*$ for the transpose and conjugate transpose of a matrix $M$. We let $I_m, O_m, J_m$ denote the $m\times m$ identity, zero and all-ones matrix. If the context is clear, then we omit the subscripts. We use boldface lowercase letters to denote vectors. We denote the empty, complete, path and cycle on $m$ vertices by $O_m$, $K_m$, and $P_m$, respectively, and the complete bipartite graph with parts of size $r$ and $s$ by $K_{r,s}$. We write $\Gamma_1\cong \Gamma_2$ to denote that $\Gamma_1$ and $\Gamma_2$ are isomorphic graphs. We also denote the automorphism group of $\Gamma$ by $\operatorname{Aut}(\Gamma)$.

Next, we review a few basic graph operations. Let $\Gamma_1$ and $\Gamma_2$ be two graphs. Let $\oplus$ and $\otimes$ denote the direct sum and Kronecker product on matrices. The \textit{union} of $\Gamma_1$ and $\Gamma_2$, denoted $\Gamma_1\sqcup\Gamma_2$, is the graph with vertex set $V(\Gamma_1) \sqcup V(\Gamma_2)$ and adjacency matrix $ A(\Gamma_1)\oplus A(\Gamma_2).$ The \textit{direct product}, \textit{Cartesian product}, and \textit{wreath product} of $\Gamma_1$ and $\Gamma_2$, denoted $\Gamma_1\times\Gamma_2$, $\Gamma_1\square\Gamma_2$, and $\Gamma_1[\Gamma_2]$, are the graphs with vertex set $V(\Gamma_1)\times V(\Gamma_2)$ and adjacency matrices
\begin{equation*}
A(\Gamma_1)\otimes A(\Gamma_2), \quad A(\Gamma_1)\otimes I+ I\otimes A(\Gamma_2),\quad \text{and}\quad  A(\Gamma_1[\Gamma_2])
=A(\Gamma_1)\otimes J+I\otimes A(\Gamma_2),
\end{equation*}
respectively, where $I$ and $J$ are the identity and all-ones square matrices of appropriate sizes.  We denote the union, direct product and Cartesian product of $k$ copies of $\Gamma$ by $\bigsqcup_k \Gamma$, $\Gamma^{\times k}$, and $\Gamma^{\square k}$, respectively. The \textit{complement} of $\Gamma$, denoted $\overline{\Gamma}$, is the graph with $V(\Gamma)=V(\overline{\Gamma})$ and
\begin{center}
$A(\overline{\Gamma})=J-I-A(\Gamma)$.
\end{center}

From \cite{BrouwerHaemers2012}, it is known that
\begin{equation}
\label{eq:graphprod}
\begin{split}
&\spec\bigl(\Gamma_1\sqcup\Gamma_2\bigr)=\spec(\Gamma_1)\cup \spec(\Gamma_2),\\
&\spec\bigl(\Gamma_1\times\Gamma_2\bigr)=\bigl\{\lambda\mu:\lambda\in \spec(\Gamma_1),\mu\in\spec(\Gamma_2)\bigr\},\\
&\spec\bigl(\Gamma_1\square\Gamma_2\bigr)=\bigl\{\lambda+\mu:\lambda\in \spec(\Gamma_1),\mu\in\spec(\Gamma_2)\bigr\}.
\end{split}
\end{equation}
If we further assume that $\Gamma_2$ is regular with valency $d$, then
\begin{equation}
\label{eq:lexprod}
\spec\bigl(\Gamma_1[\Gamma_2]\bigr)
=
\bigl\{\lambda|V(\Gamma_2)|+d:\lambda\in\spec(\Gamma_1)\bigr\}\ \bigcup\ \bigcup_{j=1}^{|V(\Gamma_1)|}\ 
\bigl\{\mu:\mu\in\spec(\Gamma_2)\setminus\{d\}\bigr\}
\end{equation}
where the union here is taken with respect to multisets.

Finally, we discuss periodicity and its relation to sedentariness. A set $S\subseteq\R$ with $|S|\geq 2$ satisfies the \textit{ratio condition} if $\frac{a-b}{c-d}\in\Q$ for all $a,b,c,d\in S$ with $c\neq d$. A non-isolated vertex $u$ of $\Gamma$ is \textit{periodic} if $|U(\tau)_{u,u}|=1$ for some $\tau>0$. The minimum such $\tau$ is called the \textit{(minimum) period} of vertex $u$ in $\Gamma$, denoted by $\rho$. It is known that vertex $u$ is periodic if and only if $\Phi_{u}$ satisfies the ratio condition \cite{Godsil2011}. We say that $\Gamma$ is a \textit{periodic} graph if every vertex in $\Gamma$ is periodic at the same time. The fact below, derived from Lemmas 9 and 12 in \cite{Monterde1}, implies that the infimum in (\ref{eq:sed}) is attained for a periodic vertex.

\begin{prop}
\label{prop:sedper}
Let $u$ be a periodic vertex with period $\rho$. The following are equivalent.
\begin{enumerate}
\item Vertex $u$ is sedentary in $\Gamma$.
\item For all $t\in[0,\rho]$, $U(t)_{u,u}\neq 0$.
\item We have $\mathscr{C}=\displaystyle\min_{t\in[0,\rho]} |U(t)_{u,u}|>0.$
\end{enumerate}
\end{prop}

We now state \cite[Corollary~13(1)]{Monterde1}, which says that vertex $u$ is sedentary whenever the $u$th diagonal entry of some $E_\lambda$ is sufficiently large. We let $\nu_2(a)$ denote the largest power of two that divides an integer $a$. We adopt the convention that $\nu_2(0)=+\infty$.

\begin{thm}
\label{thm:sed}
Let $u$ be a vertex in $\Gamma$ and fix $\lambda\in \Phi_{u}$. The following hold.
\begin{enumerate}
\item If $(E_{\lambda})_{u,u}>\frac{1}{2}$, 
then $u$ is sedentary with $|U(t)_{u,u}|\geq 2(E_{\lambda})_{u,u}-1$. 
\item Let $(E_{\lambda})_{u,u}=\frac{1}{2}$. Then $u$ is sedentary if and only if there are integers $\ell_\mu$ such that
\begin{equation}
\label{eq:pgst}
\displaystyle\sum_{\mu\in\Phi_{u}\backslash\{\lambda\}}\ell_\mu(\mu-\lambda)=0\quad  \text{and} \quad \displaystyle\sum_{\mu\in\Phi_{u}\backslash\{\lambda\}}\ell_\mu\quad \text{is odd}.
\end{equation}
Additionally, if $u$ is periodic, then $u$ is sedentary if and only if $\nu_2(\frac{\mu-\lambda}{\sqrt{\Delta}})$ are not all equal for all $\mu\in\Phi_{u}\backslash\{\lambda\}$, where $\Delta=1$ if $\Phi_{u}\subseteq\Z$ and $\Delta>1$ is a square-free integer determined from $ \Phi_{u}$ otherwise.

\end{enumerate}
\end{thm}

\medskip

\section{Twin Vertices}\label{sec:tw}
In this section, we gather known facts about twin vertices that are useful in later sections.

Let $N_{\Gamma}(u)$ denote the set of neighbors of vertex $u$ in $\Gamma$. Two vertices $u$ and $v$ are \textit{twins} if $N_{\Gamma}(u) \setminus \{u, v \} = N_{\Gamma}(v) \setminus \{u,v\}$. Further, if $u$ and $v$ are adjacent, then they are called \textit{true twins}; otherwise, they are called \textit{false twins}. The following lemma is evident. 

\begin{lem} \label{lem:joint_twins}
Let $X_1,X_2\subseteq V(\Gamma)$. If all elements of $X_1,X_2$ are pairwise true (resp., false) twins and $X_1 \cap X_2 \neq \varnothing$, then all elements of $X_1 \cup X_2$ are pairwise true (resp., false) twins. 
\end{lem}

So, if $\Gamma$ contains a vertex $u$ with a true (resp., false) twin, then we can find the largest subset $X\subseteq V(\Gamma)$ such that $u\in X$ and all elements of $X$ are pairwise true (resp., false) twins. 

A subset $X \subset V(\Gamma)$ is a \textit{twin set} if it is a maximal set in which $N_{\Gamma}(u) \setminus \{u, v \} = N_{\Gamma}(v) \setminus \{u,v\}$ for all pairs $u, v \in X$. 
In particular, if $X$ contains adjacent vertices, then $X$ is a \textit{true twin set}, and if no vertices in $X$ are adjacent, then it is a \textit{false twin set}. 
Note that $X$ is a true twin set in $\Gamma$ if and only if $X$ is a false twin set in the complement of $\Gamma.$ The following lemma will also be useful later on.

\begin{lem} \label{lem:automorphism}
Let $f\in\operatorname{Aut}(\Gamma)$. Then $X\subseteq V(\Gamma)$ is a twin set if and only if $f(X)$ is.
\end{lem}  

Next, we state a characterization of twin vertices (\cite[Lemma~2.9]{Monterde2022} and \cite[Theorem~3]{monterde2026new}).

\begin{lem}
\label{alphabeta}
Let $X$ be a twin set in $\Gamma$ with $u \in X$. The vertex $v$ belongs to $X$ if and only if
$\mathbf{e}_u-\mathbf{e}_v$ is an eigenvector for $\Gamma$ associated with the eigenvalue $\theta$, where $\theta=-1$ if $u$ and $v$ are true twins, and $\theta=0$ otherwise. Moreover, if $m=|X|$, then
\begin{equation}
\label{eq:F}
E_{\theta}=\bigg(\big(I_{m}-\frac{1}{m}J_{m}\big)\oplus O_{n-m}\bigg)+F,
\end{equation}
where $\operatorname{rank}F=m_\theta-(m-1)$, $F(\be_u-\be_v)=\zero$ for all $u,v\in X$ and $F$ is absent if $m_\theta=m-1$, where $m_\theta$ is the multiplicity of $\theta$.
\end{lem}

From the above lemma, it follows that $\theta\in\Phi_u$ whenever $u$ belongs to a twin set.

The next result establishes a dichotomous relationship between PGST and sedentariness involving twins in a graph. It follows by combining \cite[Corollary~1]{monterde2026new} and the fact that PST and PGST are equivalent for periodic vertices \cite{pal2021quantum}.

\begin{thm}
\label{thm:twdich}
Let $X$ be a twin set in $\Gamma$. Then a vertex in $X$ is either sedentary or involved in PGST with another vertex in $X$. Additionally, if a vertex in $X$ is periodic, then it is either sedentary or involved in PST with another vertex in $X$. 
\end{thm}

We now state a characterization of sedentary twin vertices first established by Monterde (see Theorems 2 and 3 in~\cite{monterde2026new}). We provide a proof for completeness.

\begin{thm}
\label{thm:twinschar}
Let $X$ be a twin set in $\Gamma$ and consider the eigenvalue $\theta$ and the matrix $F$ in Lemma~\ref{alphabeta}. Then each vertex in $X$ is sedentary if and only if one of the following holds.
\begin{enumerate}
\item $|X|\geq 3$.
\item $X=\{u,v\}$ and a $\theta$-eigenvector $\mathbf{w}$ orthogonal to $ \mathbf{e}_u-\mathbf{e}_v$ satisfies $\mathbf{w}^T\mathbf{e}_u\neq 0$.
\item $X=\{u,v\}$, $\theta$ is a simple eigenvalue of $\Gamma$, and condition (\ref{eq:pgst}) in Theorem~\ref{thm:sed} holds.
\item $X=\{u,v\}$, any $\theta$-eigenvector $\mathbf{w}$ orthogonal to $\mathbf{e}_u-\mathbf{e}_v$  satisfies $\mathbf{w}^T\mathbf{e}_u=\mathbf{w}^T\mathbf{e}_v=0$, and condition (\ref{eq:pgst}) in Theorem~\ref{thm:sed} holds.
\end{enumerate}
Moreover, if $u\in X$ and condition (1) or (2) holds, then $|U(t)_{u,u}|\geq 1-\frac{2}{|X|}+2F_{u,u}$ for all $t\in \R$, where $F$ is given in (\ref{eq:F}).
\end{thm}

\begin{proof}
Let $X$ be a twin set in $\Gamma$. If (1) or (2) holds, then each vertex in $X$ is sedentary by \cite[Theorem 3]{monterde2026new}. Now, if (3) or (4) holds, then condition (\ref{eq:pgst}) in Theorem~\ref{thm:sed} holds, and so each vertex in $X$ is sedentary by Theorem \ref{thm:sed}(2). Finally, if $X=\{u,v\}$, any $\theta$-eigenvector $\mathbf{w}$ orthogonal to $\mathbf{e}_u-\mathbf{e}_v$  satisfies $\mathbf{w}^T\mathbf{e}_u=\mathbf{w}^T\mathbf{e}_v=0$, and condition (\ref{eq:pgst}) in Theorem~\ref{thm:sed} does not hold, then $u$ and $v$ admit PGST in $\Gamma$ by \cite[Theorem 14]{kirkland2023quantum}. Thus, $u$ is not sedentary by Theorem \ref{thm:twdich}. This proves the above characterization.
\end{proof}

For the matrix $F$ in (\ref{eq:F}), it is clear from its construction in \cite{monterde2026new} that $F_{u,u}>0$ whenever condition (2) in Theorem \ref{thm:twinschar} is satisfied, while $F_{u,u}=0$ whenever $\theta$ is a simple eigenvalue of $\Gamma$, or any $\theta$-eigenvector $\mathbf{w}$ orthogonal to $ \mathbf{e}_u-\mathbf{e}_v$ satisfies $\mathbf{w}^T\mathbf{e}_u=\mathbf{w}^T\mathbf{e}_v=0$. 

From \cite[Theorem 11]{kirkland2023quantum}, we deduce the following result concerning twins with PST whose supports consist of all integers.

\begin{thm}
\label{thm:parcon}
If PST occurs between twin vertices $u$ and $v$ such that $\Phi_{u}\subseteq\Z$, then
$\nu_2(\lambda-\theta)$ is constant for all $\lambda\in\Phi_{u}\backslash\{\theta\}$, where $\theta$ is given in Lemma \ref{alphabeta}.
\end{thm}

Let $\Gamma=\Gamma_1[\Gamma_2]$ where $\Gamma_2\in\{O_m,K_m\}$. Note that all the vertices in $\{(u,j):j\in V(\Gamma_2)\}$ are twins in $\Gamma$ for each $u\in V(\Gamma_1)$. Thus, each vertex in $\Gamma$ belongs to a twin set of size at least $m$ (indeed, it will belong to a twin set that has size an integer multiple of $m$). 
Also, note that $\Gamma$ is vertex-transitive if and only if $\Gamma_1$ is vertex-transitive. To determine the support of a vertex in $\Gamma$, note that $\one$ is an eigenvector for $\Gamma_2$ associated with the eigenvalue $d=0$ if $\Gamma_2=O_m$, and the eigenvalue $d=m-1$ otherwise. If $\mathbf{x}$ is an eigenvector for the eigenvalue $\lambda$ of $\Gamma_1$, then $\mathbf{x}\otimes \one$ is an eigenvector for the eigenvalue $m\lambda+d$ of $\Gamma$. Moreover, if $\mathbf{y}$ is an eigenvector for the eigenvalue $\mu$ of $\Gamma_2$ orthogonal to $\one$, then $\mathbf{x}\otimes\mathbf{y}$ is an eigenvector for the eigenvalue $\mu$ of $\Gamma$. Therefore,
\begin{equation}
\label{eq:supplex}
\Phi_{(u,j)}(\Gamma)= 
\begin{cases}
    \{m\lambda:\lambda\in\Phi_{u}(\Gamma_1)\}\cup\{0\}, & \text{if } \Gamma_2=O_m \\
    \{m\lambda+(m-1):\lambda\in\Phi_{u}(\Gamma_1)\}\cup\{-1\},  & \text{if } \Gamma_2=K_m.
\end{cases}
\end{equation}
We use this fact to characterize the quantum state transfer properties of $\Gamma$.

\begin{thm}
\label{thm:sedwreath}
Let $\Gamma=\Gamma_1[\Gamma_2]$ where $\Gamma_2\in\{O_m,K_m\}$ and consider $\theta$ in Lemma~\ref{alphabeta}.
\begin{enumerate}
\item Vertex $(u,j)$ is involved in PGST in $\Gamma$ if and only if $\Gamma_2\in\{K_2,O_2\}$,  $\theta\notin\Phi_{u}(\Gamma_1)$,  and for all integers $\ell_\lambda$, we have
\begin{equation}
\label{eq:blowuppgst}
\displaystyle\sum_{\lambda\in\Phi_{(u,j)}(\Gamma)\backslash\{\theta\}}\ell_\lambda(\lambda-\theta)=0\quad  \text{implies that} \quad \displaystyle\sum_{\lambda\in\Phi_{(u,j)}(\Gamma)\backslash\{\theta\}}\ell_\lambda\quad \text{is even}.
\end{equation}
\item Vertex $(u,j)$ is involved in PST in $\Gamma$ if and only if $\Gamma_2\in\{K_2,O_2\}$, $\theta\notin\Phi_{u}(\Gamma_1)$, $\frac{\lambda-\theta}{\sqrt{\Delta}}$ is an integer for all $\lambda\in\Phi_{(u,j)}(\Gamma)\backslash\{\theta\}$ where $\Delta=1$ if $\Phi_{(u,j)}(\Gamma)\subseteq\Z$ and $\Delta>1$ is a square-free integer determined from  $\Phi_{(u,j)}(\Gamma)$ otherwise, and there is an integer $a$ such that $\nu_2(\frac{\lambda-\theta}{\sqrt{\Delta}})=a$ for all $\lambda\in\Phi_{(u,j)}(\Gamma)\backslash\{\theta\}.$ In this case, $\tau=\frac{\pi}{2^{a}\sqrt{\Delta}}$ is a PST time.
\end{enumerate}
Moreover, PGST or PST in $\Gamma$ occurs between pairs of vertices in the same twin set.
\end{thm}

\begin{proof}
The case when $\Gamma_2=O_m$ is known, see Theorems 3 and 4 in \cite{bhattacharjya2024quantum}. The case $\Gamma_2=K_m$ follows similarly using (\ref{eq:supplex}) and the same arguments used in Theorems 3 and 4 in \cite{bhattacharjya2024quantum}.
\end{proof}

Graphs of the form $\Gamma=\Gamma_1[\Gamma_2]$, where $\Gamma_2\in\{K_m,O_m\}$, are known in the literature as \textit{blow-up graphs}, i.e. $\Gamma$ is a graph obtained by either replacing every vertex of $\Gamma_1$ by a clique of size $m$, or replacing every vertex of $\Gamma_1$ by a coclique of size $m$; then two vertices are adjacent if they belong to cliques, or cocliques that replaced vertices adjacent in $\Gamma_1$. Quantum walks on blow-up graphs of the form $\Gamma_1[O_m]$ were investigated in \cite{bhattacharjya2024quantum,monterde2025laplacian}. The following is immediate from Theorems~\ref{thm:twdich} and~\ref{thm:sedwreath}.

\begin{cor}
\label{cor:blup}
Let $\Gamma=\Gamma_1[\Gamma_2]$ where $\Gamma_2\in\{O_m,K_m\}$ and consider $\theta$ in Lemma~\ref{alphabeta}. Vertex $(u,j)$ is sedentary in $\Gamma$ if and only if one of the following conditions holds.
\begin{enumerate}
\item $m\geq 3$.
\item $m=2$ and $\theta\in\Phi_{u}(\Gamma_1)$.
\item $m=2$, $\theta\notin\Phi_{u}(\Gamma_1)$, and condition~(\ref{eq:blowuppgst}) does not hold.
\item $m=2$, $\theta\notin\Phi_{u}(\Gamma_1)$, $\Phi_{(u,j)}(\Gamma)$ satisfies the ratio condition, and the values of $\nu_2(\frac{\lambda-\theta}{\sqrt{\Delta}})$ are not constant for all $\lambda\in\Phi_{(u,j)}(\Gamma)\backslash\{\theta\}.$
\end{enumerate}
\end{cor}

To illustrate Theorem~\ref{thm:sedwreath} and Corollary~\ref{cor:blup}, we need the next lemma.

\begin{lem}
\label{lem:cartprod}
If $u\in V(X)$ and $v\in V(Y)$, then
\begin{center}
$\Phi_{(u,v)}(X\square Y)=\{\mu+\eta:\mu\in\Phi_u(X),\eta\in\Phi_v(Y)\}$.
\end{center}
\end{lem}

\begin{proof}
If $\lambda\in\sigma(X\square Y)$, then it is known that the spectral idempotent for $\lambda$ is given by
\begin{equation*}
F_{\lambda}=\sum_{{\mu+\eta=\lambda\atop\mu\in\sigma(X),\eta\in \sigma(Y)}}E_\mu\otimes E_\eta,
\end{equation*}
where $E_\mu$ and $E_\eta$ are spectral idempotents for $\mu\in\sigma(X)$ and $\eta\in \sigma(Y)$, respectively. Thus, $(F_\lambda)_{(u,v),(u,v)}\neq 0$ if and only if at least one of the products $(E_\mu)_{u,u}(E_\eta)_{v,v}$ is nonzero. Equivalently, $\mu\in\Phi_u(X)$ and $\eta\in\Phi_v(Y)$.
\end{proof}


\begin{cor}
\label{cor:p3}
Let $n\geq 1$ and $v,w\in V(P_3^{\square n})$ such that $w=(a,\ldots,a)$, where $a$ is the degree two vertex of $P_3$. Then
\begin{center}
$\Phi_{w}(P_3^{\square n})=
\begin{cases} 
\big\{\pm 2j\sqrt{2}:j=0,1,\ldots, \frac{n}{2}\big\} & \text{if $n$ is even}\\
\big\{\pm (2j+1)\sqrt{2}:j=0,1,\ldots, \frac{n-1}{2}\big\} & \text{if $n$ is odd}.
\end{cases}$
\end{center}
Moreover, if $v\neq w$, then $\Phi_{v}(P_3^{\square n})$ is equal to $\sigma(P_3^{\square n})=\{\pm j\sqrt{2}:j=0,1,\ldots, n\}$. 
\end{cor}

\begin{proof}
Note that for all $n\geq 1$, we have
$\sigma(P_3^{\square n})=\{\pm j\sqrt{2}:j=0,1,\ldots, n\}$. 
Let $a$ and $b$ be degree two and degree one vertices of $P_3$, respectively. Using the eigenvectors of $P_3$, we get $\Phi_a(P_3)=\{\pm \sqrt{2}\}$ and $\Phi_{b}(P_3)=\{\pm \sqrt{2},0\}$. Applying Lemma \ref{lem:cartprod} to $P_3^{\square 2}=P_3\square P_3$ gives us $\Phi_{(a,a)}(P_3^{\square 2})=\{\pm 2\sqrt{2},0\}$ and the support of any other vertex in $P_3^{\square 2}$ is equal to $\sigma(P_3^{\square 2})$. Using this fact and again applying Lemma \ref{lem:cartprod} to $P_3^{\square 3}=P_3^{\square 2}\square P_3$, we get that $\Phi_{(a,a,a)}(P_3^{\square 3})=\{\pm 3\sqrt{2},\pm \sqrt{2}\}$
and the support of any other vertex in $P_3^{\square 3}$ is equal to $\sigma(P_3^{\square 3})$. One may then argue inductively to obtain the above result.
\end{proof}

We now illustrate Theorem~\ref{thm:sedwreath} and Corollary~\ref{cor:blup} in the next examples. For any graph $X$, each vertex in $X[O_m]$ belongs to a twin set of size at least $m$. Thus, if $m\geq 3$, then Corollary \ref{cor:blup}(1) implies that $X[O_m]$ is sedentary. Hence, we only deal with examples where $m=2$.

\begin{expl}
\label{ex:p3}
$K_{1,n}[O_2]$ with $n\geq 2$.\\
Let $n\geq 2$. Consider $K_{1,n}$ with central vertex $u_0$, and leaf vertices $v_1,\dots, v_n$.
Then
\begin{center}
$\Phi_{u_0}(K_{1,n})=\{\pm\sqrt{n}\}\quad \text{and}\quad \Phi_{v_j}(K_{1,n})=\{\pm\sqrt{n},0\}$.
\end{center}
Thus, for each vertex $(w, j)$ of $(K_{1,n})[O_2]$, the set 
$\Phi_{(w,j)}(K_{1,n}[O_2])=\{\pm 2\sqrt{n},0\}$ satisfies the ratio condition. In this example, we have $\theta = 0$ and $\nu_2(\frac{\lambda-\theta}{\sqrt{n}})=1$ for all $\lambda\in \Phi_{(w,j)}(K_{1,n}[O_2])\backslash\{0\}$. By 
Corollary~\ref{cor:blup}(2), each $(v_\ell, j)$ is sedentary in $K_{1,n}[O_2]$. Meanwhile, vertices $(u_0, 1)$ and $(u_0, 2)$ admit PST in $K_{1,n}[O_2]$ at time $\frac{\pi}{2\sqrt{n}}$ by Theorem \ref{thm:sedwreath}(2).
\end{expl}

\begin{figure}
	\begin{center}
		\begin{tikzpicture}[scale=.37,auto=left]
		\tikzstyle{every node}=[circle, thick, fill=white, scale=0.7]
		
		
		\node[draw,fill=cyan] (3) at (6,1.5)  {};
		\node[draw,fill=cyan] (1) at (0,1.5)  {};
		\node[draw,fill=magenta!70] (2) at (3,1.5)  {};
        \node[draw,fill=cyan] (6) at (6,-1.5)  {};
		\node[draw,fill=cyan] (4) at (0,-1.5)  {};
		\node[draw,fill=magenta!70] (5) at (3,-1.5)  {};
		
		\draw [thick, black!70] (2)--(1);
		\draw [thick, black!70] (3)--(2);
        \draw [thick, black!70] (5)--(4);
		\draw [thick, black!70] (6)--(5);
        \draw [thick, black!70] (1)--(5);
        \draw [thick, black!70] (2)--(4);
        \draw [thick, black!70] (2)--(6);
        \draw [thick, black!70] (3)--(5);
		
		\node[draw,fill=magenta!70] (11) at (15,0.85)  {};
        \node[draw,fill=magenta!70] (12) at (15,-0.85)  {};
        \node[draw,fill=cyan] (13) at (18.15,3.8)  {};
        \node[draw,fill=cyan] (14) at (18.15,2.35)  {};
        \node[draw,fill=cyan] (15) at (22.35,2.35)  {};
        \node[draw,fill=cyan] (16) at (22.35,3.8)  {};
        \node[draw,fill=magenta!70] (17) at (25.5,0.85)  {};
        \node[draw,fill=magenta!70] (18) at (25.5,-0.85)  {}; 
        \node[draw,fill=cyan] (21) at (18.15,-3.8)  {};
        \node[draw,fill=cyan] (22) at (18.15,-2.35)  {};
        \node[draw,fill=cyan] (19) at (22.35,-2.35)  {};
        \node[draw,fill=cyan] (20) at (22.35,-3.8)  {};
        
		\draw [thick, black!70] (12)--(11);
        \draw [thick, black!70] (13)--(14);
        \draw [thick, black!70] (12)--(13);
        \draw [thick, black!70] (13)--(11);
        \draw [thick, black!70] (12)--(14);
        \draw [thick, black!70] (14)--(11);
        
        \draw [thick, black!70] (15)--(16);
        \draw [thick, black!70] (17)--(18);
        \draw [thick, black!70] (15)--(17);
        \draw [thick, black!70] (15)--(18);
        \draw [thick, black!70] (17)--(16);
        \draw [thick, black!70] (18)--(16);
        
        \draw [thick, black!70] (19)--(20);
        \draw [thick, black!70] (19)--(21);
        \draw [thick, black!70] (19)--(22);
        \draw [thick, black!70] (20)--(21);
        \draw [thick, black!70] (20)--(22);
        \draw [thick, black!70] (22)--(21);

        \draw [thick, black!70] (11)--(21);
        \draw [thick, black!70] (11)--(22);
        \draw [thick, black!70] (12)--(21);
        \draw [thick, black!70] (12)--(22);

        \draw [thick, black!70] (17)--(19);
        \draw [thick, black!70] (17)--(20);
        \draw [thick, black!70] (18)--(19);
        \draw [thick, black!70] (18)--(20);
        
        \draw [thick, black!70] (13)--(16);
        \draw [thick, black!70] (13)--(15);
        \draw [thick, black!70] (14)--(16);
        \draw [thick, black!70] (14)--(15);

        \draw [thick, black!70] (11)--(17);
        \draw [thick, black!70] (11)--(18);
        \draw [thick, black!70] (12)--(17);
        \draw [thick, black!70] (12)--(18);
		
		\end{tikzpicture}
	\end{center}	
	\caption{\label{fig2}
    $P_3[O_2]$ admits PST at time $\frac{\pi}{2\sqrt{2}}$ between false twin vertices marked pink (left); $(K_2\square K_{1,2})[K_2]$ admits PGST between pairs of true twin vertices marked pink (right). All blue vertices are sedentary in both graphs.}
\end{figure}
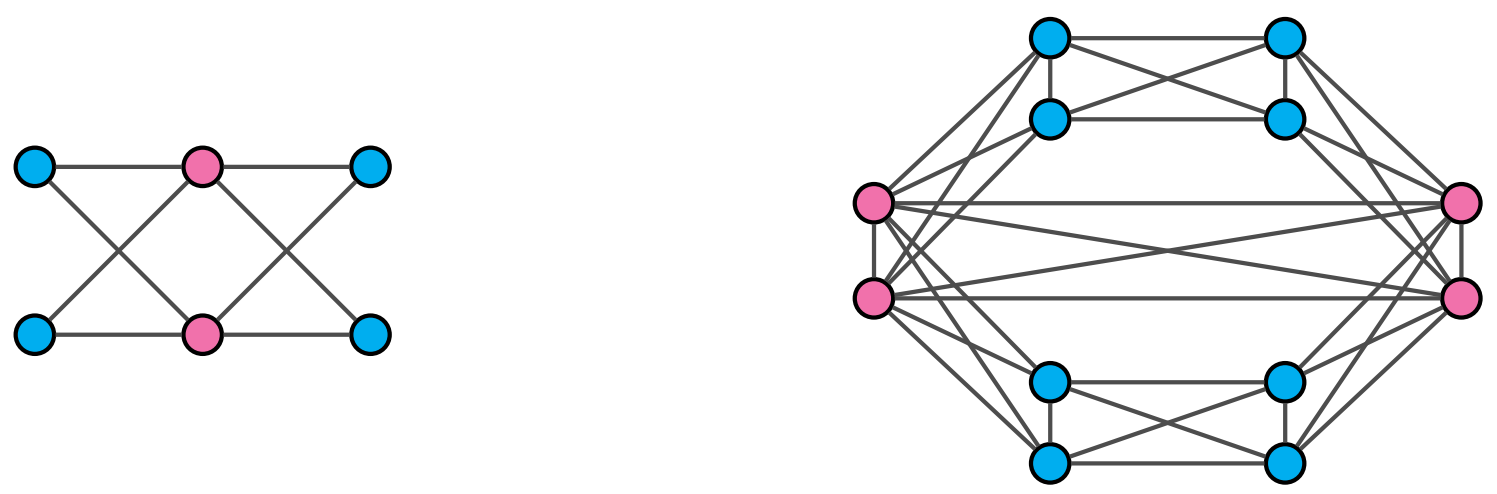

\begin{expl}$P_3^{\square n}[O_2]$  with $n\geq 1$ (see Figure \ref{fig2} for $P_3[O_2]$). \\
Let $u\in V(P_3^{\square n})$. First, suppose $n$ is odd and $u=(a,\ldots,a)$ where $a$ is the degree two vertex of $P_3$. Then by Corollary \ref{cor:p3}, we have
\begin{center}
$\Phi_{u}(P_3^{\square n})=\{\pm (2\ell+1)\sqrt{2}:\ell=0,1,\ldots, \frac{n-1}{2}\}$,
\end{center}
and so $0 \notin \Phi_{u}(P_3^{\square n})$. We have $\Phi_{(u,1)}(P_3^{\square n}[O_2])=\{\pm 2(2\ell+1)\sqrt{2}:\ell=0,1,\ldots, \frac{n-1}{2}\}\cup\{\theta=0\}$,
and so $\Phi_{(u,1)}(P_3^{\square n}[O_2])$ satisfies the ratio condition and $\nu_2(\frac{\lambda-\theta}{\sqrt{2}})$ is a constant equal to 1. Invoking Theorem \ref{thm:sedwreath}(2), we obtain PST between $(u,1)$ and its twin in $P_3^{\square n}[O_2]$ at time $\frac{\pi}{2\sqrt{2}}$.

If $n$ is even or $u\neq (a,\ldots,a)$ where $a$ is the degree two vertex of $P_3$, then $0 \in \Phi_{u}(P_3^{\square n})$ by Corollary \ref{cor:p3}. Thus, vertex $(u,j)$ in $P_3^{\square n}[O_2]$ for $j\in\{1,2\}$ is sedentary by Corollary~\ref{cor:blup}(2).
\end{expl}





\begin{expl}$(K_k\square K_{1,n})[K_2]$ for $k, n \geq 2$ (see Figure \ref{fig2} for $K_2\square K_{1,2})[K_2]$). \\
In this example we are taking the wreath product with a complete graph so $\theta = -1$. Let $w_j=(a, v_j)$ be a vertex of $K_k  \square K_{1,n}$, where $v_j$ is not the central vertex of $K_{1,n} $. Then 
$-1\in \Phi_{w_j}(K_k  \square K_{1,n} )$. By Corollary~\ref{cor:blup}(2), $(w_j, 1)$ and $(w_j, 2)$ are sedentary in $(K_k  \square K_{1,n}  )[K_2] $. Now, let $w_0=(a, u_0)$ be a vertex of $K_k \square K_{1,n}$ with $u_0$ as the central vertex of $K_{1,n}$. Then
\begin{center}
$\Phi_{w_0}(K_k  \square K_{1,n} )=\{k+\sqrt{n}-1,k-\sqrt{n}-1,\sqrt{n}-1,-\sqrt{n}-1\}$. 
\end{center}
Hence, for $j\in\{1,2\}$, $\Phi_{(w_0,j)}((K_k\square K_{1,n})[K_2])$ is equal to
\begin{center}
$\{2(k+\sqrt{n}-1)+1,2(k-\sqrt{n}-1)+1,2(\sqrt{n}-1)+1,2(-\sqrt{n}-1)+1\}\cup\{-1\}$. 
\end{center}

If $\sqrt{n} \in \Z$, then $\Phi_{(w_0,j)}((K_k\square K_{1,n})[K_2])$ satisfies the ratio condition. 
Further, if $k=\sqrt{n}$, then $-1\in \Phi_{w_0}(K_k  \square K_{1,n} )$, so $(w_0, i)$ is sedentary in 
$(K_k  \square K_{1,n} )[K_2] $ by Corollary~\ref{cor:blup}(2). 
Now, suppose $k \neq \sqrt{n}$. Here, we have $-1 \notin \Phi_{w_0}( K_k  \square K_{1,n})$, so we consider $\nu_2(\lambda-\theta)$ where $\lambda \in \Phi_{(w_0,j)}((K_k\square K_{1,n})[K_2])\backslash \{-1\}$. This is constant for all $\lambda$ if and only if $\nu_2(k)>\nu_2(\sqrt{n})$.  Thus, if $\nu_2(k)>\nu_2(\sqrt{n})$, then 
Theorem~\ref{thm:sedwreath}(2) implies that vertex $(w_0, 1)$ in $(K_k  \square K_{1,n} )[K_2]$ admits PST with its twin at time $\frac{\pi}{2g}$ where $g=\operatorname{gcd}(k+\sqrt{n},k-\sqrt{n},\sqrt{n})$.  If $\nu_2(k)\leq \nu_2(\sqrt{n})$, then by Corollary~\ref{cor:blup}(4), vertex $(w_0, i)$ is sedentary in $(K_k  \square K_{1,n})[K_2]$ for $i\in\{1,2\}$.

If $\sqrt{n} \notin \Z$, then $-1\notin \Phi_{w_0}(K_k  \square K_{1,n})$, so we consider condition~(\ref{eq:blowuppgst}).
Set 
$$
0= \ell_1(k+\sqrt{n})+\ell_2(k-\sqrt{n})+\ell_3\sqrt{n}-\ell_4\sqrt{n}=(\ell_1+\ell_2)k+(\ell_1-\ell_2+\ell_3-\ell_4)\sqrt{n},
$$ 
where each $\ell_j$ is an integer. As $\sqrt{n}\notin\Z$, we get $\ell_1+\ell_2=\ell_1-\ell_2+\ell_3-\ell_4=0$, and thus $\ell_1+\ell_2+\ell_3+\ell_4=2(\ell_1+\ell_3)$ is always even. By Theorem~\ref{thm:sedwreath}(1), PGST occurs between $(w_0,1)$ and its twin in $(K_k \square K_{1,n} )[K_2]$.

\end{expl}

\section{Walk-regular graphs}
\label{sec:wr}

A graph $\Gamma$ is \textit{$k$-regular} if every vertex has degree $k$, this is equivalent to $\mathbf{1}$ being an $k$-eigenvector for $\Gamma$. If $\Gamma$ is regular, the eigenvalues of $\overline{\Gamma}$ can be determined from the eigenvalues of $\Gamma$.
Let $U_\Gamma(t)$ and $U_{\overline{\Gamma}}(t)$ denote the transition matrices of $\Gamma$ and $\overline{\Gamma}$, respectively. The following result is due to Godsil~\cite[Lemma~2.1]{godsil2021sedentary}.

\begin{lem}
\label{reg}
If $\Gamma$ is a regular graph on $n$ vertices, then for all $u,v\in V(\Gamma)$,
\begin{equation*}
\big|U_{\overline{\Gamma}}(t)_{u,v}-e^{-it}U_\Gamma(-t)_{u,v}\big|\leq 2/n\quad \text{for all $t\in\R$}.
\end{equation*}  
\end{lem}

\begin{thm}
\label{regcomp}
Let $\Gamma$ be a regular graph on $n$ vertices with vertex $u$ that is $C$-sedentary. If $C>\frac{2}{n}$, then vertex $u$ is $(C-\frac{2}{n})$-sedentary in $\overline{\Gamma}$.
\end{thm}

\begin{proof}
By assumption, $|U_\Gamma(-t)_{u,u}|\geq C$ for all $t\in\R$. Since $\Gamma$ is regular, Lemma~\ref{reg} yields
\begin{equation*}
|U_{\overline{\Gamma}}(t)_{u,u}|\geq \big|U_\Gamma(-t)_{u,u}\big|- 2/n\geq C-2/n>0\quad \text{for all $t\in\R$}.\qedhere
\end{equation*}
\end{proof}

A regular graph $\Gamma$ is \textit{walk-regular} if each nonnegative power of $A(\Gamma)$ has a constant diagonal. Examples of walk-regular graphs include vertex-transitive graphs, Cayley graphs and graphs in association schemes (such as distance-regular graphs).

We will now show that for walk-regular graphs, the diagonal entries of $U(t)$ can be written as a function in $e^{\rm{i}t}$.

\begin{prop}
\label{prop:walk-regular}
Let $\Gamma$ be a walk-regular graph on $n$ vertices. For any vertex $u \in V(\Gamma)$, 
\begin{equation*}
U(t)_{u,u}=\frac{1}{n} \sum_{j=1}^n e^{\rm{i} \lambda_j t}\quad \text{for all $t\in\R$}.
\end{equation*}
\end{prop}

\begin{proof}
By definition, we have 
\begin{equation*}
U(t)=\sum_{k=0}^{\infty} \frac{(\rm{i} t)^k}{k!} A^k\quad \text{for all $t\in\R$}. 
\end{equation*}
As $\Gamma$ is walk-regular, $A^k$ has constant diagonal for all $k \geq 0.$ Thus, $U(t)$ has constant diagonal for all $t\in\R$, and so $U(t)_{u,u}$ is equal to $\frac{1}{n}$ times the trace of $U(t)$. If $\lambda_1\geq \cdots\geq \lambda_n$ are the eigenvalues of $\Gamma$, then the eigenvalues of $U(t)$ are $e^{\rm{i} \lambda_j t}$, and so the conclusion follows.
\end{proof}

\begin{rem}
\label{rem:suppwr}
If $m_\lambda$ is the multiplicity of $\lambda\in \sigma(\Gamma)$, then for any vertex $u\in V(\Gamma)$, Proposition \ref{prop:walk-regular} yields
\begin{equation*}
U(t)_{u,u}=\frac{1}{n} \sum_{\lambda\in \sigma(\Gamma)}m_\lambda e^{\rm{i} \lambda t}\quad \text{for all $t\in\R$}.
\end{equation*}

Consequently, if $\Gamma$ is walk-regular, since each $E_\lambda$ is a polynomial in $A(\Gamma)$, each $E_\lambda$ has constant diagonal. As the trace of $E_\lambda$ is equal to the multiplicity $m_\lambda$ of $\lambda$, we get that $(E_\lambda)_{u,u}=\frac{m_\lambda}{n} > 0$ for all $\lambda\in\sigma(\Gamma)$, thus $\Phi_{u}=\sigma(\Gamma)$ for all $u\in V(\Gamma)$.
\end{rem}

From Proposition \ref{prop:walk-regular}, the following is immediate.

\begin{prop}
\label{prop:sedwrg}
A walk-regular graph is sedentary if and only if one of its vertices is.
\end{prop}

Let $\Gamma_1\vee \Gamma_2$ denote the join of $\Gamma_1$ and $\Gamma_2$, this is the graph obtained by taking one copy of each of $\Gamma_1$ and $\Gamma_2$, and adding all edges between $\Gamma_1$ and $\Gamma_2$. We remark there are regular graphs, that are not walk-regular and are sedentary. One example is $\Gamma=(\bigsqcup_p K_3)\vee (\bigsqcup_q K_4)$ where $4q-3p=1$ and $p,q>1$. Here, each vertex of $\Gamma$ belongs to a twin set of size at least three, which is sedentary by Theorem~\ref{thm:twinschar}(1). There are also regular but not walk-regular graphs that admit a subset $W\subseteq V(\Gamma)$ such that each vertex in $W$ is sedentary but each vertex in $V(\Gamma)\backslash W$ is not sedentary. One example is $\Gamma=(\bigsqcup_{2p} K_6)\vee (\bigsqcup_{2q} K_2)$ with $3p-q=1$, where the vertices that belong to $\bigsqcup_{2p} K_6$ are all sedentary in $\Gamma$ by Theorem~\ref{thm:twinschar}(1), while those belonging to $\bigsqcup_{2q} K_2$ admit PST by \cite[Theorem~2]{Alvir2016} (see also~\cite[Theorem 6]{kirkland2026quantum}). 

By remark \ref{rem:suppwr}, we may restate Theorem~\ref{thm:sed} in the context of walk-regular graphs.

\begin{cor}
\label{cor:sed}
Let $\Gamma$ be a  walk-regular graph and $\lambda\in\sigma(\Gamma)$. 
\begin{enumerate}
\item If $m_\lambda>\frac{n}{2}$, 
then $\Gamma$ is sedentary with $|U(t)_{u,u}|\geq \frac{2m_\lambda}{n}-1$. 
\item If $m_\lambda=\frac{n}{2}$, then $\Gamma$ is sedentary if and only if condition (2) of Theorem~\ref{thm:sed} holds.
\end{enumerate}
\end{cor}

The complement of a walk-regular graph is walk-regular. The disjoint union of $k$ copies of a walk-regular graph is also walk-regular. So, if $\Gamma$ is a walk-regular graph, then so are $\bigsqcup_k \Gamma$ and $\overline{\bigsqcup_k \Gamma}$ for all $k\geq 1$. This observation combined with Theorem~\ref{regcomp} yields the following.

\begin{cor}
\label{cor:regcomp1}
Let $\Gamma$ be a $C$-sedentary walk-regular graph on $n$ vertices. If $C>\frac{2}{n}$, then $\overline{\bigsqcup_k \Gamma}$ is a $(C-\frac{2}{n})$-sedentary walk-regular graph for all $k\geq 1$.
\end{cor}

\medskip

\section{Integral graphs}\label{sec:integ}

We say that $\Gamma$ is an \emph{integral graph} if all of its eigenvalues are integers. By the ratio condition, an integral graph is periodic. Moreover, since the valency of a regular graph is an integer and the set of eigenvalues are algebraic integers, the ratio condition implies that a regular graph is periodic if and only if it is integral \cite{coutinho2014quantum}. 

\subsection{The necklace polynomial of an integral graph}

Inspired by \Cref{prop:walk-regular} and recent works~\cite{hyde2022cyclotomic, minh2026cyclotomic}, we introduce the notion of the necklace polynomial of an integral graph. We first recall the general framework for the definition of this polynomial.  Let $R=\Z[\Z^{\circ}]$ be the algebra associated with the monoid $(\Z, \times)$. Elements of $R$ are of the form $\sum_{m \in \Z} a_m [m]$ where $m \in \Z$ and $a_m \in \Z.$ Addition and multiplication are defined naturally:
\[ 
\sum_{m} a_m[m] + \sum_{m} b_m [m] = \sum_{m} (a_m+b_m)[m], \quad  [m][n]=[mn],
\]
and extended linearly.
If $\Gamma$ is an integral graph with eigenvalues $\lambda_1 \geq  \ldots \geq \lambda_n$, then we  define
\[ \alpha_\Gamma := \sum_{j=1}^n [\lambda_j]\in R.\]

\begin{expl}
Since $\sigma(K_n)$ is $\{n-1^{(1)}, -1^{(n-1)} \}$, we have $\alpha_{K_n} = [n-1]+(n-1)[-1].$
\end{expl}

The following proposition is immediate from the spectra of graph products in (\ref{eq:graphprod}). 

\begin{prop}
Let $\Gamma_1, \Gamma_2$ be integral graphs. Then $\Gamma_1 \times \Gamma_2$, $\Gamma_1 \sqcup \Gamma_2$, and $\Gamma_1 \square \Gamma_2$ are integral graphs. Furthermore, $\alpha_{\Gamma_1 \sqcup \Gamma_2} = \alpha_{\Gamma_1} + \alpha_{\Gamma_2}$ and $\alpha_{\Gamma_1 \times \Gamma_2} = \alpha_{\Gamma_1} \alpha_{\Gamma_2}.$
\end{prop}

Given $\alpha \in \Z[\Z^{\circ}]$ with $\alpha = \sum_{m} a_m [m]$, we define the \textit{necklace function} associated with $\alpha$ as 
\[ \alpha \cdot x = \sum_{m} a_m x^m. \]
Let $m_0$ be the minimal value of $m$ appearing in $\alpha.$ The \textit{necklace polynomial} associated with $\alpha$ is given by
\[ \widehat{\alpha \cdot x} = \sum_{m} a_m x^{m-m_0}.\]

Let $z_0\in\C$ such that $|z_0|=1$. We write $(\alpha \cdot x)(z_0)$ for the evaluation of $\alpha \cdot x$ at $x=z_0.$ This is a well-defined value. Note that $z_0$ is a root of $(\alpha \cdot x)$ if and only if it is a root of $\widehat{\alpha \cdot x}.$
\begin{definition}
Let $\Gamma$ be an integral graph with eigenvalues $\lambda_1\geq \ldots\geq\lambda_n$. The \emph{necklace polynomial $N_{\Gamma}(x)$ of $\Gamma$} is the necklace polynomial of $\alpha_\Gamma.$ More concretely,
    \[ 
      N_\Gamma(x) = \sum_{j=1}^n x^{\lambda_j- \lambda_n} = x^{-\lambda_n} \sum_{j=1}^n x^{\lambda_j}  = x^{-\lambda_n} \alpha_\Gamma (x).
       \]
\end{definition}

\begin{expl}
    If $\Gamma=K_n$, then $\N_\Gamma(x) = x^n +(n-1).$
\end{expl}

\begin{prop}
    Let $\Gamma$ be a connected, integral, walk-regular graph. Then $N_\Gamma(x)$ is a monic polynomial of degree $\lambda_1-\lambda_n.$ Furthermore $N_\Gamma (1)=|V(\Gamma)|.$ 
\end{prop}

Let $P(x)\in\C[x]$ be a polynomial of degree $m$. The \textit{reciprocal polynomial} of $P(x)$ is the polynomial $P^{\star}(x):= x^m P(1/x)$. Note that the coefficients of $P^{\star}(x)$ are the coefficients of $P(x)$ in reverse order. Thus, $P(x)$ and $P^{\star}(x)$ have the same degree whenever $P(x)$ has nonzero constant term. Also, observe that $\lambda\in\C \setminus \{0\}$ is a root of $P(x)$ if and only if $1/\lambda$ is a root of $P^{\star}(x)$. We say that $P(x)$ is \textit{reciprocal} if $P(x)=P^{\star}(x)$.

\begin{lem}
\label{lem:biprec}
$N_{\Gamma}(x)$ is reciprocal if and only if $\Gamma$ is a bipartite graph. 
\end{lem}

\begin{proof}
First, observe that
    \[ N_{\Gamma}^{\star}(x) = x^{\lambda_1 - \lambda_n} \sum_{j=1}^n (1/x)^{\lambda_j-\lambda_n} = \sum_{j=1}^n x^{\lambda_1 - \lambda_j}.\]
Thus, $N_{\Gamma}(x)$ is reciprocal if and only if $\lambda_1-\lambda_j = \lambda_{n+1-j} - \lambda_n$ for all $1 \leq j \leq n$. That is, $\lambda_1 + \lambda_{n} = \lambda_{j} + \lambda_{n+1-j}$ for all $1 \leq j \leq n$. Since $\sum_{j=1}^n \lambda_j=0$, the preceding equations hold if and only if $\lambda_{j} =- \lambda_{n+1-j}$ for each $1 \leq j \leq n$. Equivalently, $\Gamma$ is bipartite.
\end{proof}

The following result establishes a connection between the sedentariness of an integral walk-regular graph and the arithmetic properties of its
necklace polynomial.

\begin{thm}
\label{thm:rootsedentariness}
 Let $\Gamma$ be an integral walk-regular graph. The following are equivalent.
 \begin{enumerate}
     \item $\Gamma$ is sedentary.
     \item $N_\Gamma(x)$ has no  roots on the unit circle.
     \item $N_\Gamma^{\star}(x)$ has no  roots on the unit circle.
 \end{enumerate}
\end{thm}

\begin{proof}
By assumption, $\Gamma$ is periodic. Moreover, Proposition \ref{prop:walk-regular} implies that for all $u\in V(\Gamma)$, 
     \[ U(t)_{u,u} =\frac{1}{n} \sum_{j=1}^n e^{\rm{i} \lambda_j t} =\frac{1}{n} \sum_{j=1}^n  z^{\lambda_j} = \frac{1}{n} z^{\lambda_n} N_{\Gamma}(z)\]
where $z = \exp(\mathrm{i} t)$. 
Suppose $\Gamma$ is sedentary. As $\Gamma$ is periodic, Proposition \ref{prop:sedper}(3) gives us $\mathscr{C}=\min_{t \in [0, 2 \pi]}|U(t)_{u,u}|>0$. Thus, $|U(t)_{u,u}|=\frac{|N_{\Gamma}(z)|}{n}\geq \mathscr{C}$ for all $t \in [0, 2 \pi]$, and so $|N_{\Gamma}(z)| \geq \mathscr{C}n>0$ for all $|z|=1.$ Hence, $N_{\Gamma}(x)$ has no root on the unit circle. Conversely, suppose $N_{\Gamma}(x)$ has no root on the unit circle. 
If $U(t_0)_{u,u}=0$ for some $t_0 \in [0, 2 \pi]$, then $\exp(it_0)$ is a root of $\alpha \cdot x$, and hence, $\exp(it_0)$ is a root of $N_{\Gamma}(x),$ a contradiction. Thus, $U(t)_{u,u}\neq 0$ for all $t\in[0,\rho]$. Since $\Gamma$ is periodic, $\Gamma$ is sedentary by Proposition \ref{prop:sedper}(1). Thus, (1) and (2) are equivalent. The equivalence of (2) and (3) follows from the fact that a unit complex number $\alpha$ is a root of $N_{\Gamma}(x)$ if and only if $1/\alpha$ is a root of $N_{\Gamma}^{\star}(x)$.
\end{proof}

\begin{rem}
By Lemma~\ref{lem:biprec}, if $\Gamma$ is an integral nonbipartite walk-regular graph, then we have two polynomials, $N_\Gamma(x)$ and $N_\Gamma^{\star}(x)$, at our disposal to test sedentariness in $\Gamma$.
\end{rem}

\begin{cor}
\label{cor:sedgcd}
Let $\Gamma$ be an integral walk-regular graph. If $\operatorname{gcd}(N_\Gamma(x),N_\Gamma^{\star}(x))=1$, then $\Gamma$ is sedentary. 
\end{cor}

\begin{proof}
We will prove the contrapositive. Suppose $\Gamma$ is not sedentary. Then $N_\Gamma(x)$ has a root $\lambda$ on the unit circle by Theorem~\ref{thm:rootsedentariness}. Since $N_\Gamma(x)$ is a polynomial with integer coefficients, it follows that all its conjugates, including $\frac{1}{\lambda}$, are also roots of $N_\Gamma(x)$. Since $\frac{1}{\lambda}$ is a root of $N_{\Gamma}(x)$, it follows that $N_{\Gamma}(x)$ and $N_{\Gamma}^{\star}(x)$ share a common non-constant factor.
\end{proof}

\begin{rem}
\label{rem:1}
The converse of Corollary~\ref{cor:sedgcd} need not be true. We provide a counterexample in the family of unitary Cayley graphs introduced in \cite{klotz2007some} and generalized in \cite{unitary, nguyen2025gcd}. Let $R$ be a finite commutative ring. Let $S$ be a subset of $R$ such that $S=-S$. We say that $\Gamma(R,S)$ is a \textit{gcd-graph} if $S$ is stable under  multiplication by units in $R.$ When $S=R^{\times}$, the graph is called a \textit{unitary Cayley graph} and we write $G_{R}=\Gamma(R,R^{\times}).$ Gcd-graphs possess elegant properties, such as the explicit description of their spectra via Ramanujan sums. In particular, it is known that they are integral graphs. These properties lie at the cross of algebra, number theory, and graph theory which makes it an ideal candidate for our investigation.

Let $\Gamma = G_{12}  = \Cay(\Z_{12}, \{1,5,7,11\})$ (the \emph{unitary Cayley graph} on $\Z_{12}$). Then we have 
\[ N_{\Gamma}(x)=x^8 + 2x^6 + 6x^4 + 2x^2 + 1.\]

By the triangle inequality, one can verify that $N_{\Gamma}(x)$ has no root on the unit circle, and so $\Gamma$ is sedentary by Theorem~\ref{thm:rootsedentariness}. This polynomial is reciprocal, so $\operatorname{gcd}(N_\Gamma(x),N_\Gamma^{\star}(x))\neq 1$. 
\end{rem}

\subsection{An algebraic lemma}

In this section, we describe some special cases where we can use tools from algebraic number theory to show that $N_{\Gamma}(x)$ has no root on the unit circle. 
\begin{lem} \label{lem:root-poly}
Let $P(x) \in \Z[x]$ be a monic polynomial of degree $m$. Suppose that $P(x)$ can be written as $P(x)=x^m + q Q(x)$ where $q \geq 2$ is a positive integer  and $Q(x) \in \Z[x]$ is a polynomial of degree less than $m.$ Then, $P(x)$ has no root on the unit circle. 
\end{lem}
We provide two proofs for \Cref{lem:root-poly}. They are based on the same idea. The first proof is shorter but the second proof could be used to give an explicit bound on sedentariness. 

\begin{proof}[First proof]
By way of contradiction, suppose $\theta \in \C$ such that $|\theta|=1$ and $P(\theta)=0.$ Let $F$ be the splitting field of $P$ and $\OO_{F}$ be its ring of algebraic integers. By assumption, $\theta \in \OO_F.$ Since $P(x)$ has real coefficients, $\bar{\theta}$ is also a root of $P(x).$ Thus, $\bar{\theta} \in \OO_F.$ Since $|\theta|=1$, we have $\theta \bar{\theta}=1$, and so $\theta$ is a unit in $\OO_F.$ Now because $\theta^m = -q Q(\theta)$ and $\theta$ is a unit in $\OO_F$, it follows that $-q$ is a unit in $\OO_F$. Since $q \geq 2$, this is a contradiction. 
\end{proof}

\begin{proof}[Second proof]
    Let $P^{\star}(x) := x^m P(1/x)$ be the reciprocal polynomial of $P.$ We first claim that 
    \[ \gcd(P, P^{*})=1 \] 
    and the result will follow from the proof for Corollary~\ref{cor:sedgcd}.
    Suppose to the contrary that $P$ and $P^{*}$ have a nonconstant common factor, say $F(x) \in \Q[x].$ Since $P$ is monic, by Gauss's lemma on primitive polynomials, we may assume that $F(x)$ is a monic polynomial in $\Z[x].$ Let us write $P= FG_1, P^{\star}=F G_2$ where $G_1, G_2 \in \Z[x]$ both have degree less than $m$.  Let $p$ be a prime divisor of $q$. Taking reduction modulo $p$, we have the following equalities in $\F_p[x]$ 
    \[ x^m = \bar{F} \bar{G_1}, 1 = \bar{F} \bar{G_2}.\]
    Since $F$ is monic, the second equation implies that $\bar{F}$ is a constant while the first equation implies that $\bar{F}=x^k$ for some $1 \leq k \leq m.$ This is a contradiction. 
\end{proof}

\begin{thm} \label{thm:main-multiple-eig}
Let $\Gamma$ be a connected, integral, walk-regular graph with valency $d$. If the multiplicity of each $\lambda\in\sigma(\Gamma)\backslash\{d\}$ is a multiple of some integer $q\geq 2$, then $\Gamma$ is sedentary.
\end{thm}

\begin{proof}
By our assumption, $N_{\Gamma}(x)= x^m + q Q(x)$ where $\deg(Q)<m.$ By \cref{lem:root-poly}, we know that $N_{\Gamma}(x)$ has no root on the unit circle. Thus, $\Gamma$ is sedentary by \cref{thm:rootsedentariness}. 
\end{proof}

\begin{cor}
\label{cor:dirprod}
Let $\Gamma = K_{n_1} \times \cdots \times K_{n_r}$. If $\gcd(n_1-1, \ldots, n_r-1)>1$, then $\Gamma$ is sedentary. In particular, $K_n^{\times k}$ is sedentary for all $n \ge 3$ and $k \geq 1$.
\end{cor}
\begin{proof}
The valency of $\Gamma$ is given by $d=\prod_{i=1}^r(n_i-1).$  For a subset $I\subseteq\{1,\ldots,r\}$, we have that $(-1)^{|I|}\prod_{i\notin I}(n_i-1)$ is an eigenvalue of $\Gamma$ with multiplicity $\prod_{i\in I}(n_i-1)$.
Let $q=\gcd(n_1-1,\ldots,n_r-1)>1$. Then, every nonempty subset contributes a multiplicity divisible by $q$. Furthermore, if several subsets produce the same eigenvalue, their contributions are added. Consequently, every eigenvalue other than $d$ has multiplicity
divisible by $q$. We then see that Theorem~\ref{thm:main-multiple-eig} yields the result. 
\end{proof}

Let $\Gamma=K_{n_1} \times \cdots \times K_{n_r}$ where $n_j=2$ for some $j$. Then $\gcd(n_1-1, \ldots, n_r-1)=1$ and so the conclusion in Corollary~\ref{cor:dirprod} does not apply. Nevertheless, since $\Gamma$ in this case is a nonsingular bipartite graph, \cite[Theorem~12]{meagher2026sedentary} implies that $\Gamma$ is not sedentary.

\subsection{Explicit bound for sedentariness}
\Cref{thm:main-multiple-eig} implies that if the necklace polynomial of $\Gamma$ has a certain form, then $\Gamma$ is sedentary. In this section, we find an explicit constant $C$ such that $\Gamma$ is $C$-sedentary. To do this, we will use the theory of resultants.

Let $f= a_nx^n + a_{n-1}x^{n-1}+\cdots+a_0 \in \C[x]$  be a polynomial. The \textit{height} of $f$ is defined as
\[ H(f) = \max_{0 \leq i \leq n } |a_i|.  \]

\begin{lem} \label{lem:bound-solution}
    Let $M \in M_N(\C)$ be an invertible matrix. Let $\x =(x_1, \ldots, x_N)^T \in \C^N$ be a solution to the equation $M\x = \bv$, where $\bv$ is a unit vector. Let $B$ be a positive number such that $|M_{ij}| \leq B$ for all $1 \leq i, j \leq N.$ Then, for each $1 \leq k \leq N$, we have 
 \[ |x_k| \leq \frac{B^{N-1} N^{(N-1)/2}}{|\det(M)|}. \]
 Furthermore, if $M$ has integer coefficients, then $|x_k| \leq B^{N-1} N^{(N-1)/2}$ for each $1\ \leq k \leq N.$

\end{lem}
\begin{proof}
Let $\bc_1, \ldots, \bc_N$ be the column vectors of $M.$ By our assumption, $ \|\bc_i\| \leq B \sqrt{N}$
for each $1 \leq i \leq N$.
    By Cramer's rule, $x_k =  \dfrac{\det(M_k)}{\det(M)}$ where $M_k$ is the matrix obtained by replacing the $k$-th column of $M$ by $\bv.$ Applying Hadamard's inequality, we have 
    \[ |\det(M_k)| \leq \|\bv\|\cdot  \prod_{j \neq k} \|\bc_j\| \leq  B^{N-1} N^{\frac{N-1}{2}}. \]
    Therefore, 
     \[ |x_k| \leq \frac{B^{N-1} N^{(N-1)/2}}{|\det(M)|}. \]
     If $M$ has integer coefficients then $|\det(M)| \geq 1$ and the inequality follows. 
\end{proof}

For polynomials $f,g\in\C[x]$, we let $\deg(f)$ denote the degree of $f$ and $\Res(f,g)$ denote the resultant of $f$ and $g.$ 

\begin{prop} \label{prop:bound-coefficients}
Let $f, g \in \C[x]$ be polynomials such that $\gcd(f,g)=1$. Let $B\in\R$ such that $B \geq \max \{H(f), H(g)\}.$ Then, there exist $p, q \in \C[x]$ such that:
    \begin{enumerate}
        \item $pf + qg =1.$
        \item $\deg(p) \leq deg(g)-1$ and $\deg(q) \leq deg(f)-1.$
        \item If $N= deg(f)+deg(g)$, then 
        \[ H(p), H(q) \leq \dfrac{B^{N-1}N^{(N-1)/2}}{|\Res(f,g)|}. \] 
        In particular, if $f, g \in \Z[x]$ then we have $ H(p), H(q) \leq B^{N-1} N^{(N-1)/2}.$
    \end{enumerate}
\end{prop}

\begin{proof}
    The coefficients of $p,q$ are solutions of a linear system of the form $S\x = \be_N$ where $S$ is the Sylvester matrix of $f$ and $g$, and $\be_N$ is a standard vector. So the  
the existence of $p, q$, with the properties (1) and (2), follows from elimination theory, see~\cite[Proposition~9, Chapter 5]{cox2007ideals}. Furthermore, by definition, $\Res(f,g) = \det(S)$, so the upper bound on $H(p)$ and $H(q)$ then follows from \cref{lem:bound-solution}.
\end{proof}

\begin{rem}
The bound in \cref{prop:bound-coefficients} is not optimal, since each column of $S$ is obtained from the coefficients of $f$ and $g$, and zeros. Thus, the bound on $\|\bc\|$, where $\bc$ is a column of $S$, can be sharpened. However, for our purposes, the above bound is sufficient.
\end{rem}


\begin{lem} 
\label{lem:bound-C}
Let $B>0$ and $\mathcal{P}_{n,B} \subset \Z[x]$ be a set of polynomials such that for all $P \in \mathcal{P}_{n,B}$ the following hold:
\begin{enumerate}
    \item $\deg(P)=n$ and $H(P) \leq B$ and $P(0) \neq 0;$ and
    \item  $\gcd(P, P^{\star})=1 $, where $P^{*}(x)= x^n P(1/x)$ is the reciprocal polynomial of $P$.
\end{enumerate}
Then there exists a constant $C>0$, depending only on $n$ and $B$, such that for all $z \in \C$ with $|z|=1$, we have  $|P(z)| \geq C$ for all $P \in \mathcal{P}_{n, B}.$ In fact, we may choose $C=\frac{1}{B^{2n-1} (2n)^{(2n+1)/2}}.$
\end{lem}

\begin{proof}
   Let $P \in \mathcal{P}_{n,B}$. Since $\gcd(P, P^{\star})=1$, by \cref{prop:bound-coefficients}, we can find $p, q$ such that $\deg(p), \deg(q) \leq n-1$ such that $p P + qP^{\star}=1$ and $H(p), H(q) \leq B^{2n-1} (2n)^{(2n-1)/2}.$  

   Let $z \in \mathbb{C}$ with $|z|=1$. Since $P$ has real coefficients, $P(1/z)=P(\overline{z})=\overline{P(z)}$. Thus, $|P^\star(z)|=|z^nP(1/z)|=|P(z)|$. By the triangle inequality, $|p(z)| \leq n H(p) \leq  n  B^{2n-1} (2n)^{(2n-1)/2}$.
    Similarly, $|q(z)| \leq n  H(q) \leq  n  B^{2n-1} (2n)^{(2n-1)/2}.$ Since $1 = P(z) p(z)+ P^{\star}(z) q(z)$, we get
\[ 1 \leq |P(z)\|p(z)| + |P^{\star}(z)| |q(z)| = |P(z)| (|p(z)|+|q(z)|).\]
by the triangle inequality. Consequently, we have 
\[ |P(z)| \geq \frac{1}{|p(z)|+|q(z)|} \geq \frac{1}{2n  B^{2n-1} (2n)^{(2n-1)/2}}=\frac{1}{B^{2n-1} (2n)^{(2n+1)/2}}.\qedhere \]
\end{proof}

We are now ready to prove the main result in this subsection.

\begin{thm}
\label{thm:C}
Let $\Gamma$ be an integral walk-regular graph on $n$ vertices such that $N_{\Gamma}(x)$ has the form $x^m + q Q(x)$ for some integer $q \geq 2$ and polynomial $Q(x)\in\mathbb{Z}[x]$ such that $\deg(Q)<m.$ Then $\Gamma$ is $C$-sedentary with $C = \frac{1}{n^{2n} (2n)^{(2n+1)/2}}$.
\end{thm}

\begin{proof}
As $\deg(Q(x))<m$, we have $\deg(N_\Gamma)=m$ and we may choose $m=\lambda_1-\lambda_n$. By \cite[Proposition 1.2]{gregory2001spread}, $m \leq n.$ We know that $\gcd(N_{\Gamma}(x), N_{\Gamma}^{\star}(x))=1$ by the second proof of \Cref{lem:root-poly}. Since the multiplicity of each eigenvalue of $\Gamma$ is bounded above by $n$, we also know that $H(N_\Gamma) \leq n. $ Thus, for all $t \in \R$ and for every vertex $u$ of $\Gamma$, \Cref{lem:bound-C} yields
    \[ |U(t)_{u,u}| = \frac{1}{n} |N_\Gamma(z)| \geq \frac{1}{n}\cdot \frac{1}{n^{2m-1} (2m)^{(2m+1)/2}} = \frac{1}{n^{2m} (2m)^{(2m+1)/2}} \geq \frac{1}{n^{2n} (2n)^{(2n+1)/2}}.\qedhere \]
\end{proof}

\medskip

\section{Cayley graphs}\label{sec:cayley}

We now apply our results to Cayley graphs. In particular, our main goal in this section is to show that an integral Cayley graph on a group of odd order is sedentary.

\subsection{Spectra of Cayley graphs}

Let $G$ be a finite group and $S\subseteq G$. The \textit{Cayley graph} $\Gamma(G,S)$ has vertex set $G$ and $g,h\in G$ are adjacent if and only if $g^{-1} h\in S$. Throughout, we assume $S$ is a nonempty inverse-closed subset and $1 \not \in S$, so that $\Gamma(G,S)$ is a simple undirected graph with at least one edge. Note that $\Gamma(G,S)$ is connected if and only if $S$ generates $G$.
Let $\rho: G \to GL_n(\C)$ be an irreducible unitary representation of $G$ of degree $n= d_{\rho}=\dim(\rho)$. Then $\rho$ extends to a ring map $\rho: \C[G] \to M_n(\C).$ For $S \subset G$, define 
$$
     \rho_S := \sum_{s \in S} \rho(s).
$$
We have the following theorem about the spectrum of a Cayley graph. 

\begin{thm}\cite[Theorem 3]{MR626813}\label{thm:speccayley}
The eigenvalues of a Cayley graph $\Gamma(G,S)$ are the union of the multisets of eigenvalues of $\rho_S$, each repeated $d_\rho$ times, where $\rho$ runs over one representative from each equivalence class of irreducible unitary representations of $G$. Further, if the multiplicity of $\lambda$ as an eigenvalue of $\rho_S$ is $m(\lambda,\rho)$, then its multiplicity as an eigenvalue of $\Gamma(G,S)$ is equal to $\sum_\rho d_\rho m(\lambda,\rho)$.
\end{thm}

We let $\rho_0$ be the trivial representation of $G$. Note that 
$\trivialS = |S|$ is the valency of $\Gamma(G,S).$ 

\begin{thm} \label{thm:odd-order}
Let $\Gamma=\Gamma(G,S)$ be a Cayley graph, where $G$ is a finite group of odd order. Then the multiplicity of each 
$\lambda\in\sigma(\Gamma)\backslash\{ \trivialS \}$ is an even number 
while that of $\trivialS$ is odd. 
\end{thm}

\begin{proof}
Let $G'$ be the group generated by $S.$ Then $\Gamma(G,S)$ is a disjoint union of $[G:G']$ copies of $\Gamma(G',S).$ As $[G:G']$ is odd, it suffices to prove the result for the connected graph $\Gamma(G',S)$.

Suppose $\Gamma(G,S)$ is connected so that $\trivialS$ has multiplicity one. Let $\lambda$ be another eigenvalue of $\Gamma(G,S)$ so that $\lambda<\trivialS$. By Theorem \ref{thm:speccayley}, $\lambda$ is an eigenvalue of $\rho_S$ where $\rho: G \to \GL_n(\C)$ is a non-trivial irreducible unitary representation of $G$. Let $\bar{\rho}$ be the conjugate representation of $\rho$; namely $\bar{\rho}: G \to \GL_n(\C)$ is defined by $\bar{\rho}(g)= \overline{\rho(g)}.$ Since $\rho$ is unitary, $\rho(s^{-1}) = \rho(s)^{-1} = \rho(s)^{*}= \overline{\rho(s)^T}.$ Therefore, $\bar{\rho}(s^{-1})= \rho(s)^T.$
Since $|G|$ is odd,  $\rho \neq \bar{\rho}$ by the Frobenius-Schur index theorem, see \cite[Proposition 39]{serre1977linear}. Using the fact that $S=S^{-1}$, we have 
$$
   \bar{\rho}_S = \sum_{s \in S} \bar{\rho}(s) 
                = \sum_{s \in S} \bar{\rho}(s^{-1})
                = \sum_{s \in S} \rho(s)^T
                =\rho_S^T.
$$
Now, because $\rho_S$ and $\rho_S^T$ have the same characteristic polynomial, this shows that $\lambda$ also appears in $\bar{\rho}_S$. Since $\rho \neq \bar{\rho}$ and $\rho_S$ and $\bar{\rho}_S$ have the same characteristic polynomial, the non-trivial irreducible representations occur in distinct conjugate pairs contributing the same multiplicity of $\lambda$. Hence the total multiplicity of $\lambda$ is even.
\end{proof}

Since a Cayley graph is necessarily walk-regular,  we obtain our main result in this section.

\begin{thm}
\label{thm:oddsed}
Let $G$ be a finite group of odd order. If $\Gamma(G,S)$ is integral, then it is sedentary. 
\end{thm}

\begin{proof}
Let $G'=\langle S\rangle$. Then $\Gamma(G,S)$ is a disjoint union of $[G:G']$ copies of the connected Cayley graph $\Gamma(G',S)$. Since $G'$ has odd order, \Cref{thm:odd-order} implies that every eigenvalue of $\Gamma(G',S)$ other than the degree has even multiplicity. So by \Cref{thm:main-multiple-eig}, $\Gamma(G',S)$ is sedentary. Thus, every connected component of $\Gamma(G,S)$ is sedentary, and so is $\Gamma(G,S)$.
\end{proof}

\begin{rem}
Later, we will show that if $G$ is a finite group of even order, then we can always find $S$ such that $\Gamma(G,S)$ is not sedentary (see \Cref{cor:even-not-sedentary}). 
\end{rem}

\subsection{Sedentariness in gcd-graphs}

In this subsection, we study sedentariness of graphs defined over a finite ring $R$ of characteristic $n$. In this case, the underlying vertex set is $(R,+)$. As observed in various previous works (see \cite{klotz2007some, nguyen2025gcd, so2006integral}) what makes working with graphs defined over $R$ more interesting is the interplay between the multiplicative structure and the additive structure of $R$. In particular, we have a quite elegant theory for rings in which $\widehat{R}$ is canonically isomorphic to $R$ (see \cite{honold2001characterization, nguyen2026u}).

\begin{definition}
A finite ring $R$ of characteristic $n$ is called a \textit{symmetric Frobenius ring} if there exists a generating character $\chi: (R, +) \to \C^{\times}$ satisfying (i) $\chi$ is symmetric; namely $\chi(ab)=\chi(ba)$ for all $a, b \in R$; and (ii) the map $R \to \widehat{R}=\Hom(R,\C^{\times})$ given by $r \mapsto \chi_r$, where $\chi_r(a)=\chi(ra)$, is an isomorphism. 
\end{definition}
Examples of symmetric Frobenius rings include finite quotients of a Dedekind domain, finite fields, matrix rings, finite chain rings, and their products (see \cite{honold2001characterization, nguyen2025integral}). 
 
 Let $\Gamma(R,S)$ be a Cayley graph where the underlying vertex set is $(R,+)$. When $R$ is a finite symmetric Frobenius ring, the spectrum of $\Gamma(R,S)$ is canonically parametrized by $R$. More precisely, the spectrum of $\Gamma(R,S)$ is the multiset $\{\lambda_r \}_{r \in R}$ where
 \[ \lambda_r = \sum_{s \in S} \chi_r(s)=\sum_{s \in S} \chi(rs).\]

 We will use this explicit description to find  instances where the conditions of \cref{thm:main-multiple-eig} apply. To do so, we recall the following definition which was introduced in \cite{nguyen2026u}. 
 
\begin{definition}
\label{def:6.6}
Let $U$ be a subgroup of $R^{\times}.$ A Cayley graph $\Gamma(R,S)$ is \textit{$U$-unitary} if $S$ is stable under the left and right actions of $U$, i.e.~$USU=S.$ If $\Gamma(R,S)$ is $U$-unitary and $U=R^{\times}$, then $\Gamma(R,S)$ is called a \textit{gcd-graph}. 
\end{definition}

Since $R$ is a finite $\Z_n$-algebra, there is an induced group homomorphism $(\Z_n)^{\times} \to R^{\times}.$ In what follows, we will identify $(\Z_n)^{\times}$ with its image in $R.$ It is known that the Cayley graph $\Gamma(R,S)$ is integral if and only if $\Gamma(R,S)$ is $(\Z_n)^{\times}$-unitary; namely $S$ is stable under the multiplication by elements of  $(\Z_n)^{\times}$ (see \cite{godsil2025integral, nguyen2025supercharacters}). 


\begin{lem} \label{lem:equal-eigenvalues}
Let $\Gamma(R,S)$ be a $U$-unitary Cayley graph. If $a=ubv$ for some $u,v \in U$, then $\lambda_a = \lambda_b.$ 
\end{lem}

\begin{proof}
We have
    \[ \lambda_a = \sum_{s \in S} \chi_a(s) = \sum_{s \in S} \chi(as) =\sum_{s \in S} \chi(ubvs)=\sum_{s \in S}\chi(bvsu) = \sum_{t \in vSu} \chi(bt) = \lambda_b. \]
    Here, the last equality follows from the fact that $vSu=S.$
\end{proof}
In particular, the eigenvalue $\lambda_a$ is equal for every $a$ belonging to a $U$-double orbit in $R$, so the multiplicity of any nontrivial eigenvalue is a sum of the cardinalities of such orbits. Combining \cref{thm:main-multiple-eig} and \cref{lem:equal-eigenvalues}, we have the following result. 

\begin{prop} \label{prop:double-coset}
Let $R$ be a finite symmetric Frobenius ring and $\Gamma(R,S)$ be a connected integral $U$-unitary Cayley graph where $U\leq R^\times$. Suppose $q\geq 2$ is an integer such that the cardinality of every double orbit $UaU$ with $a\neq 0$ is divisible by $q$. Then $\Gamma(R,S)$ is sedentary.
\end{prop}

We now discuss a few cases where we can apply \Cref{lem:equal-eigenvalues} and \Cref{thm:main-multiple-eig}. 

\begin{prop} \label{prop:F_q-algebra}
  Let $q=p^f\geq 3$ be a prime power and $R$ be a finite symmetric Frobenius $\F_q$-algebra. If $\Gamma(R,S)$ is a connected $\F_q^{\times}$-unitary Cayley graph, then the following hold. 
  \begin{enumerate}
      \item $\Gamma(R,S)$ is integral. 
      \item The order of every orbit in $\F_q^{\times} \backslash R^{*} / \F_q^{\times}$ is exactly $q-1$. Here $R^{*} = R \setminus \{0\}.$
      \item $\Gamma(R,S)$ is sedentary. 
  \end{enumerate}
\end{prop}

\begin{proof}
Since $S$ is stable under the action of $\F_q^{\times}$, it is also stable under the action of $\F_p^{\times}$ since $\F_p^{\times} \subset \F_q^{\times}.$ The integrality of $\Gamma(R,S)$ follows from  (1) follows from \cite[Theorem~1.1]{nguyen2025integral}. Meanwhile, (2) follows from the fact that the orbit of an element $a \in R^{*}$ in the double-quotient $\F_q^{\times} \backslash R^{*} / \F_q^{\times}$ is precisely $\{ca \mid c \in \F_q^{\times}\}$ (since $\F_q$ is in the center of $R$, we only need to consider one-sided action). Finally, (3) follows from (1), (2) and \Cref{thm:main-multiple-eig}.
\end{proof}

\begin{expl}
    Let $R=M_2(\F_4)$ and $\Gamma:=\Gamma(M_2(\F_4), \GL_2(\F_4))$ be a unitary Cayley graph. The characteristic polynomial of $\Gamma$ is $(x - 180) (x + 12)^{75}  (x - 4)^{180}$. The eigenvalues $-12$ and $4$ have multiplicity divisible by $3=4-1$ as predicted by \Cref{prop:F_q-algebra}. We conclude that the unitary Cayley graph over $R$ is sedentary. 
\end{expl}

Here is another situation where we can apply \Cref{thm:main-multiple-eig}.

\begin{prop}
Let $R$ be a finite commutative Frobenius ring; namely $R= \prod_{i=1}^d R_i$ where each $R_i$ is a local commutative Frobenius ring. Let $\kappa_i$ be the size of the residue field of $R_i.$ Suppose there is an integer $q \geq 2$ such that $\kappa_i \equiv 1 \pmod{q}$ for each $1 \leq i \leq d.$ If  $\Gamma(R,S)$ is a connected $R^{\times}$-unitary Cayley graph, then the following statements hold. 
    \begin{enumerate}
      \item $\Gamma(R,S)$ is integral. 
      \item The order of every orbit in $R^{\times} \backslash R^{*} / R^{\times}$ is divisible by $q$. Here $R^{*} = R \setminus \{0\}.$
      \item $\Gamma(R,S)$ is sedentary. 
  \end{enumerate}
\end{prop}

\begin{proof}
Note that (1) is clear, while (3) follows from (1) and (2). So we only need to prove (2). By \cite[Corollary~2.7]{nguyen2025gcd}, we know that for each $a \in R^{*}$, the orbit containing $a$ is precisely $(R/\rm{Ann}_{R}(a))^{\times}$ where $\rm{Ann}_{R}(a)=\{x \in R \mid ax=0\}.$ Since $a \neq 0$, $\rm Ann_{R}(a)$ is a proper ideal of $R$; namely, it has the form $\prod_{i=1}^d I_i$ where $I_i$ is an ideal in $R_i$ and there exists an index $k$ such that $I_k \neq R_k.$ In this case, $|(R/\rm{Ann}_{R}(a))^{\times}|$ is a multiple of $|(R_k/I_k)^{\times}|.$ By definition, the residue field of $R_k/I_k$ is also the residue field of $R_k.$ Therefore $|(R_k/I_k)^{\times}|$ is a multiple of $\kappa_k-1$ and hence a multiple of $q$. Thus $|(R/\rm{Ann}_{R}(a))^{\times}|$ is a multiple of $q$ for each $a \in R^{*}.$ 
\end{proof}


\begin{expl}
Let $n=pq$ where $p,q \equiv 1 \pmod{3}$ are distinct primes and let $S=\{a \in \Z_n \mid \gcd(a,n) \in \{ p, q\} \}.$ The characteristic polynomial of the gcd-graph 
$\Gamma(\Z_n,S)\cong K_p \square K_q$ is given by $(x - (p+q-2) )  (x - (p-2))^{q-1} (x - (q-2))^{p-1}  (x + 2)^{(p-1)(q-1)}.$  We see that all eigenvalues except $p+q-2$ have multiplicity divisible by $3$. More generally, if $n = \prod_{i=1}^r p_i^{a_i}$ where $p_i$ is a prime such that $p_i \equiv 1 \pmod{3}$. Then any gcd-graph over $\Z_n$ will be sedentary.
\end{expl}

\medskip

\section{Twins in Cayley graphs}\label{sec:twcayley}

Our goal in this section is twofold. First, we characterize the existence of twins in Cayley graphs and vertex-transitive graphs, and second, we use this characterization to completely determine sedentariness in vertex-transitive graphs with twins. Throughout this section, $G$ denotes a finite group.

\subsection{Structure theorem}

A subset $X\subseteq V(\Gamma)$ is a \textit{homogeneous set} if every vertex in $V(\Gamma)\backslash X$ is adjacent to either all or none of the vertices in $X$ \cite[Definition 2.11]{chudnovsky2024prime}. Thus, $X$ is a true (resp., false) twin set if and only if $X$ is a maximal homogeneous set in $\Gamma$ and the subgraph induced by $X$ is complete (resp., empty). 
In \cite[Theorem~3.4]{chudnovsky2024prime}, the authors provide a characterization of Cayley graphs that possess non-trivial homogeneous sets. In this section, we study the same problem for the existence of a twin set, which is somewhat easier since a twin set is necessarily a homogeneous set. Since $X$ is a true twin set in $\Gamma$ if and only if $X$ is a false twin set in $\overline{\Gamma}$, it is sufficient to focus on true twin sets in $\Gamma$. The following fact is known.

\begin{lem}
\label{lem:twaut}
For each $g \in G$, the left multiplication map $\ell_g: \Gamma(G,S) \to \Gamma(G,S)$ defined by $x \mapsto gx$ is a graph automorphism.
\end{lem}


We are now ready to prove one of the main results about twins in Cayley graphs. 

\begin{prop} \label{prop:existence_of_twins}
If $X$ is a true twin set in $\Gamma(G,S)$ and $1 \in X$, then $X$ is a subgroup of $G$. The same is true if $X$ is a false twin set. 
\end{prop}

\begin{proof}
Let $h \in X$. By Lemmas \ref{lem:automorphism} and \ref{lem:twaut}, $hX$ is a twin set. Since $1\in X$, we get that $h \in hX \cap X$, and so all elements in $hX \cup X$ are pairwise true twins by Lemma~\ref{lem:joint_twins}. By the maximality of $X$ and $hX$, we have $hX=X.$ Thus $X$ is closed under multiplication. Now, let $h \in X.$ Then, there exists $n>0$ such that $h^n=1.$ Therefore, $h^{-1} = h^{n-1} \in X$, and so $X$ is closed under inverses. This shows that $X$ is a subgroup of $G.$
\end{proof}

\begin{lem} 
\label{lem:conditions_for_twins}
    Let $H$ be a subgroup of $G$ and $\Gamma(G,S)$ be a Cayley graph. Let $S^{*} = S \cup \{1\}.$ Let $g_1H, \ldots, g_rH$ be the set of all left cosets of $H$ in $G.$ The following are equivalent.
    \begin{enumerate}
        \item Pairwise, elements of $H$ are true twins.
        \item $S^{*}$ is a union of cosets of $H$ in $G$ and $H \subset S^{*}.$
    \end{enumerate}
    Moreover, if either of the above conditions holds, then $\Gamma(G,S)\cong\Gamma_1 [K_{m}]$ where $m=|H|$ and the vertices of $\Gamma_1$ are the left cosets of $H$, and $g_i H$ and $g_j H$ are adjacent if and only if $g_i^{-1} g_j \in S$ (this is independent of the choice of $g_i$ for each coset as explained in \Cref{rem:indepenence}).  
\end{lem}

\begin{proof}
Assume every pair of elements in $H$ are twins. Since the subgraph induced by $H$ is complete, $(1,h) \in E(\Gamma)$ for all $ h\in H \setminus \{1\}.$ Since $h \in S$ by definition, we get $H \subset S^{*}.$ Next, we claim that if $g \in S$, then $gH \subset S^{*}.$ If $g \in H$, then we are done, so suppose that $g \not \in H$. Since $S$ is closed under inversion, $g^{-1} \in S$ which implies $(g^{-1}, 1) \in E(\Gamma).$ From the fact that $H$ is a twin set, we know that $(g^{-1}, h) \in E(\Gamma)$ for all $h \in H$. Therefore, $gh \in S$ for all $h \in H.$ In other words $gH \subset S \subset S^{*}.$

Conversely, suppose that $H \subset S^{*}$ and $S^{*}$ is a union of left cosets of $H$ in $G.$  We claim that the elements of $H$ are twins. First, we remark that since $H \subset S^{*}$, the induced subgraph on $H$ is the complete graph. Next, we show that $H$ is a homogeneous set. In fact, let $g \in G \setminus H.$  Then $(g,h) \in E(\Gamma)$ for some $h \in H$ if and only if $g^{-1}h \in S.$ Since $S^{*}$ is a union of left-cosets, we know this happens if and only if $g^{-1} H \subset S$, which implies $g$ is adjacent to every element in $H$. Thus, $H$ is a homogeneous set in $\Gamma.$

The last statement follows from the fact that each coset $g_iH$ is a set of twins in $\Gamma.$
\end{proof}

\begin{rem} \label{rem:indepenence}

Suppose $S^*=S\cup\{1\}$ is a union of left cosets of $H$. Then $S^*$ is a union of double cosets of $H$. Indeed, if $g\in S^*$, then $gH\subseteq S^*$. Since $(S^*)^{-1}=S^*$ (because $S$ is inverse closed), we also have $Hg\subseteq S^*$. Hence $HgH\subseteq S^*$ for every $g\in S^*$. 

We claim that the adjacency relation on the left cosets described in Lemma ~\ref{lem:conditions_for_twins}  is well-defined. Indeed, if $g_i'=g_i h_i$ and $g_j'=g_j h_j$ with $h_i,h_j\in H$, then $(g_i')^{-1}g_j'=h_i^{-1}g_i^{-1}g_jh_j$. Since $S^*$ is a union of double cosets of $H$, we have $g_i^{-1}g_j\in S^*$ if and only if $(g_i')^{-1}g_j'\in S^*$. 
\end{rem}

\begin{rem}

Further note that $H$ is twin set if is it maximal, so $H$ is a twin set if and only if $H$ is a maximal subgroup with $H \subset S^*$ so that $S^*$ is a union of cosets of $H$.
\end{rem}

\begin{rem}
\label{rem:Gamma1}
If $H$ is a normal subgroup of $G$, then $\Gamma_1$ in Lemma  ~\ref{lem:conditions_for_twins} is a Cayley graph $\Gamma_1=\Gamma(G/H,S)$ over the quotient group $G/H$, where $S$ is the set of all cosets $g_iH$ 
such that $g_i H \subset S^{\star}$ and $g_i\notin H$.

\end{rem}
    
If we replace ``true twins" by ``false twins'' in \Cref{lem:conditions_for_twins}, then $\Gamma$ is isomorphic to $\Gamma_1[O_m]$, where $O_m$ denotes the empty graph on $m=|H|$ vertices.

\begin{lem}
\label{lem:conditions_for_twins2}
With the same premise in Proposition \ref{lem:conditions_for_twins}, the following are equivalent.
\begin{enumerate}
\item Pairwise, elements of $H$ are false twins.
\item $H\cap S=\varnothing$ and $S$ is a union of left cosets of $H$. 

\end{enumerate}
Moreover, if one of the above conditions holds, then $\Gamma\cong \Gamma_1[O_m]$, where $m=|H|$ and the vertices of $\Gamma_1$ are the left cosets of $H$, and $g_iH$ and $g_jH$ are adjacent if and only if $g_i^{-1}g_j\in S$.
\end{lem}

\begin{rem}
\label{rem:verttrans}
In Lemmas~\ref{lem:conditions_for_twins} and \ref{lem:conditions_for_twins2},  $\Gamma_1$ is also vertex transitive. To see this, note that the group $G$ acts on the vertices of $\Gamma_1$ (although this action is not faithful) and induces a transitive action on $\Gamma_1$.

In general if $\Gamma = \Gamma_1 [\Gamma_2]$ where $\Gamma_2 = K_m$ or $O_m$ forms a twin set, then $\Gamma$ is vertex transitive if and only if $\Gamma_1$ is. To see this, suppose $X_i$ and $X_j$ are two vertices of $\Gamma_1$, and let $x_i\in X_i$ and $x_j\in X_j$ in $\Gamma$. Since $\Gamma$ is vertex transitive, there exists $f\in\operatorname{Aut}(\Gamma)$ such that $f(x_i)=x_j.$ By Lemma~\ref{lem:automorphism} and maximality of twin sets, we have $ f(X_i)=X_j.$ Moreover, $f$ permutes $X_1,\ldots,X_r$ and preserves adjacency between them. Thus, $f$ induces an automorphism $\overline{f}$ of $\Gamma_1$ such that $ \overline{f}(X_i)=X_j.$ Therefore $\Gamma_1$ is vertex transitive. Clearly if $\Gamma_1$ is vertex transitive, then $\Gamma$ is as well.
\end{rem}

The proof of \Cref{lem:conditions_for_twins} can be extended to vertex-transitive graphs.  

\begin{thm}
\label{thm:conditions_for_twins_vert-trans}
Let $\Gamma$ be a graph with a twin set $X$ of size $m$. Then $\Gamma$ is vertex transitive if and only if there is a vertex-transitive graph $\Gamma_1$ such that $\Gamma\cong \Gamma_1[\Gamma_2]$, where $\Gamma_2=K_m$
if $X$ is a true twin set, and $\Gamma_2=O_m$ otherwise. In particular, $V(\Gamma)$ can be partitioned into sets of twins of cardinality $m$.
\end{thm}

\begin{proof}
The converse is clear, so we only prove the forward direction. We start with the case when $X$ is a true twin set in $\Gamma$. By Lemmas \ref{lem:automorphism} and \ref{lem:joint_twins}, we get that for every $f\in \text{Aut}(\Gamma)$, either $
f(X)=X$ or $f(X)\cap X=\varnothing.$ 
Now fix $x\in X$. Since $\Gamma$ is vertex-transitive, for every $v\in V(\Gamma)$ there exists $f\in\operatorname{Aut}(\Gamma)$ such that $f(x)=v$. Therefore
\[
V(\Gamma)=\bigsqcup_{f\in\operatorname{Aut}(\Gamma)}f(X).
\]
Let $f, g \in \text{Aut}(\Gamma)$ such that $ f(X)\cap g(X)\neq\varnothing.$
Then $X\cap f^{-1}g(X)\neq \varnothing$
and therefore, we get $ f^{-1}g(X)=X.$
Thus $f(X)=g(X)$, so the distinct automorphic images of $X$ form a partition
\[
V(\Gamma)=X_1\cup \cdots\cup X_r,
\]
where each $X_i$ is a true twin set and $ |X_i|=|X|=m.$ Define a graph $\Gamma_1$ with vertex set $V(\Gamma_1)=\{X_1,\ldots,X_r\}$ and
two vertices $X_i$ and $X_j$ are adjacent in $\Gamma_1$ if there exists an edge in $\Gamma$ between a vertex of $X_i$ and a vertex of $X_j$.
This is well-defined since each $X_i$ is a twin set. Since the subgraph induced by $X_i$ is $K_m$, we have $\Gamma\cong\Gamma_1[K_m].$ If $X$ is instead a false twin set, the same argument gives a partition into sets $X_i$ of false twins of cardinality $m$. In this case, the subgraph induced by each $X_i$ is the empty graph $O_m$, and hence $\Gamma\cong\Gamma_1[O_m].$ In both cases, the same argument used in Remark~\ref{rem:verttrans} yields the fact that $\Gamma_1$ is vertex-transitive.
\end{proof}

\subsection{Applications to sedentariness}

We now characterize sedentariness in vertex-transitive graphs with twins. 

\begin{thm}
\label{thm:caytwsed}
Let $\Gamma$ be a vertex-transitive graph with a twin set $X$ such that $u,v\in X$. Consider $\theta$ in Lemma~\ref{alphabeta}. Then $\Gamma$ is sedentary if and only if $\Gamma\cong \Gamma_1[\Gamma_2]$ for some graph $\Gamma_1$ with $n$ vertices, $\Gamma_2\in\{K_m,O_m\}$ where $m=|X|\geq 2$, and one of the following holds:
\begin{enumerate}
\item $m\geq 3$;
\item $m=2$ and $\theta\in\spec(\Gamma_1)$; 
\item $m=2$, $\theta\notin\spec(\Gamma_1)$, and there are integers $\ell_\lambda$ such that
\begin{equation*}
\displaystyle\sum_{\lambda\in\sigma(\Gamma)\backslash\{\theta\}}\ell_\lambda(\lambda-\theta)=0\quad  \text{and} \quad \displaystyle\sum_{\lambda\in\sigma(\Gamma)\backslash\{\theta\}}\ell_\lambda\quad \text{is odd};
\end{equation*}
\item $m=2$, $\theta\notin\spec(\Gamma_1)$, $\Gamma$ is integral, and the values of $\nu_2(\lambda-\theta)$ are not constant for all $\lambda\in\sigma(\Gamma)\backslash\{\theta\}.$
\end{enumerate}
Moreover, if condition (1) or (2) holds, then $\Gamma$ is $C$-sedentary, where $C=1-\frac{2}{m}+\frac{2(m_\theta-n(m-1))}{mn}$ and $m_\theta$ is the multiplicity of $\theta$ as an eigenvalue of $\Gamma$.
\end{thm}

\begin{proof}
Our assumption combined with  \Cref{thm:conditions_for_twins_vert-trans} implies that $\Gamma\cong \Gamma_1[\Gamma_2]$ for some vertex-transitive graph $\Gamma_1$ with $n$ vertices and $\Gamma_2\in\{K_m,O_m\}$, where $m=|X|\geq 2$. Because of vertex transitivity, Remark \ref{rem:suppwr} implies that the support of each vertex in $\Gamma$ and $\Gamma_1$ equals $\sigma(\Gamma)$ and $\sigma(\Gamma_1)$, respectively. Applying Corollary \ref{cor:blup} yields the result.
\end{proof}

A similar argument used in Theorem \ref{thm:caytwsed} combined with Theorem \ref{thm:sedwreath} yields a characterization of PST and PGST in vertex-transitive graphs with twins.

\begin{cor}
\label{cor:pstvertextranstw}
Let $\Gamma$ be a vertex-transitive graph with twins
and consider $\theta$ in Lemma~\ref{alphabeta}.
\begin{enumerate}
\item PGST occurs in $\Gamma$ if and only if $\Gamma\cong \Gamma_1[\Gamma_2]$ for some vertex-transitive graph $\Gamma_1$ such that $\theta\notin\spec(\Gamma_1)$, $\Gamma_2\in\{K_2,O_2\}$, and for all integers $\ell_\lambda$, we have
\begin{equation*}
\displaystyle\sum_{\lambda\in\sigma(\Gamma)\backslash\{\theta\}}\ell_\lambda(\lambda-\theta)=0\quad  \text{implies that} \quad \displaystyle\sum_{\lambda\in\sigma(\Gamma)\backslash\{\theta\}}\ell_\lambda\quad \text{is even}.
\end{equation*}
\item PST occurs in $\Gamma$ if and only if $\Gamma\cong \Gamma_1[\Gamma_2]$ for some vertex-transitive graph $\Gamma_1$ such that $\theta\notin\sigma(\Gamma_1)$, $\Gamma_2\in\{K_2,O_2\}$, 
$\Gamma$ is integral, and there is an integer $a$ such that $\nu_2(\lambda-\theta)=a$ for all $\lambda\in\sigma(\Gamma)\backslash\{\theta\}.$ In this case, $\tau=\frac{\pi}{2^{a}}$ is a PST time.
\end{enumerate}
Moreover, if PGST or PST occurs in $\Gamma$, then it occurs between every pair of twins in $\Gamma$.
\end{cor}

We now illustrate Theorem~\ref{thm:caytwsed} and Corollary~\ref{cor:pstvertextranstw} in the following examples. 

\begin{expl}
\label{ex:qd}
Let $Q_d$ be the hypercube of dimension $d\geq 1$. It is known that
\begin{equation*}
\sigma(Q_d)=\{d-2j:j=0,1,\ldots,d\},
\end{equation*}
and so $Q_d$ is integral, and hence periodic. Moreover $0\notin \sigma(Q_d)$ if and only if $d$ is odd, while $-1\notin \sigma(Q_d)$ if and only if $d$ is even.
Making use of (\ref{eq:lexprod}), we get
\begin{center}
$\sigma(Q_d[O_2])=\{2(d-2j):j=0,1,\ldots,d\}\cup \{\theta=0\}$
\end{center}
and
\begin{center}
$\sigma(Q_d[K_2])=\{2(d-2j)+1:j=0,1,\ldots,d\}\cup \{\theta=-1\}$.
\end{center}
For $Q_d[O_2]$, we have $\nu_2(\lambda-\theta)=\nu_2(2(d-2j)-0)$ is a constant $a$ for all $j\in\{0,1,\ldots,d\}$ if and only if $d$ is odd, in which case $a=1$. For $Q_d[K_2]$, we have $\nu_2(\lambda-\theta)=\nu_2(2(d-2j+1)$ is a constant $a$ for all $j\in\{0,1,\ldots,d\}$ if and only if $d$ is even, in which case $a=1$. Applying Corollary \ref{cor:pstvertextranstw}(2), PST occurs in $Q_d[O_2]$ (resp., $Q_d[K_2]$) if and only if $d$ is odd (resp., even), in which case, it occurs between every pair of twins at time $\frac{\pi}{2}$. Moreover, by Theorem \ref{thm:caytwsed}, $Q_d[O_m]$ (resp., $Q_d[K_m]$) is sedentary if and only if $m\geq 3$ or $d$ is even (resp., odd).
\end{expl}

\begin{expl}
As $\sigma(K_n[O_2])=\{2(n-1),-2\}\cup\{\theta=0\}$, the same argument used in Example~\ref{ex:qd} yields PST in $K_n[O_2]$ if and only if $n$ is even, in which case, it occurs between every pair of twins at  $\frac{\pi}{2}$. Also, $K_n[O_m]$ is sedentary if and only if $m\geq 3$ or $n$ is odd.
\end{expl}

\begin{expl}
Let $C_n$ be the cycle on $n$ vertices. It is known that $C_n[O_2]$ does not admit PGST \cite[Theorem 7]{bhattacharjya2024quantum}. Thus each vertex in $C_n[O_m]$ is sedentary for all $m\geq 2$ by Theorem ~\ref{thm:twdich}.
\end{expl}




We now complement Theorem~\ref{thm:oddsed} by showing that if $|G|$ is even, then there exists a subset $S \subset G$ such that $\Gamma(G,S)$ is connected and not sedentary. 

\begin{cor}
\label{cor:even-not-sedentary}
If $|G|$ is even, then there exists a subset $S$ such that $\Gamma(G,S)$ is connected and not sedentary. 
\end{cor}

\begin{proof}
First, suppose $|G|=4k$ for some integer $k.$ Let $z\in G$ such that $z^2=1$, $H=\{1, z\}$ and consider $\Gamma = \Gamma(G,S)$ where $S= G \setminus H$. By \Cref{lem:conditions_for_twins2}, $\Gamma \cong K_{2k}[O_2]$. Since $\Gamma$ is a cocktail party graph, it is not sedentary by~\cite[Corollary~4]{monterde2026new}. Next, let $|G|=2k$ where $k$ is odd. Let $z\in G$ such that $z^2=1$. We first show that $G$ has a normal subgroup of index $2$. The action on $G$ by itself via left multiplication induces a group homomorphism $\ell: G \to S_{2k}.$ For each pair $(g,zg)$, $z$ acts as a transposition. The sign of $\ell(z)$ is therefore equal to $(-1)^k =-1.$ This shows that the composition map $\ell: G \to S_{2k} \to \pm{1}$ is surjective. Let $K$ be the kernel of this map. Then, $K$ is a normal subgroup of index $2$ in $G.$ Let $S= (K \setminus \{1\}) \cup \{z\}$. Then, $S$ is inverse closed and any two elements in $K$ are adjacent. Similarly, any two elements in $zK$ are adjacent. Furthermore, each element in $K$ is adjacent to exactly one element in $zK$. Thus, $\Gamma(G,S)\cong K_2 \Box K_{k}$, which is not sedentary. In both cases, $\Gamma$ is connected. 
\end{proof}

\begin{rem}
Without the connectedness requirement, we may simply take $S= \{z\}$ where $z^2=1.$ In this case, $\Gamma(G,S)$ is a disjoint union of copies of $K_2$, and hence not sedentary. 
\end{rem}



\medskip

\section{Modular properties of eigenvalues}
\label{sec:modprop}

In Theorem~\ref{thm:parcon}, we showed that for PST to occur between twins $u$ and $v$ in a graph $\Gamma$, it must be that $\nu_2(\lambda-\theta)$ is equal for all $\lambda\in\Phi_{u}\backslash\{\theta\}$. In this penultimate section, we investigate whether this modular property holds in general for integral graphs and Cayley graphs. 

\subsection{Integral graphs}

We now examine whether the above-mentioned modular property holds for integral graphs.

\begin{lem} \label{lem:congruence_condition_join}
Let $\Gamma$ be an integral graph with eigenvalues $\lambda_1,\ldots, \lambda_n$. If $n=|V(\Gamma)|$ is odd, then the following statements hold.
\begin{enumerate}
    \item There is no  $a\in\Z$ such that $\nu_2(\lambda_j)=a$ for all $1 \leq j \leq n.$
    \item If some $\lambda_i$ is odd and $\lambda_j\neq -1$ for all $1\leq j\leq n$, then there is no $a \in \Z$ such that $\nu_2(\lambda_j+1)=a$ for all $1 \leq j\leq n.$
    \item If $\lambda_i=0$ has even multiplicity, then there is no $a\in\Z$ such that $\nu_2(\lambda_j)=a$ for all $\lambda_j\neq 0$. 
    \item If $\lambda_i=-1$ has even multiplicity and for some $\ell\neq i$, we have $\lambda_\ell$ is odd and $\lambda_\ell\neq -1$, then there is no $a \in \Z$ such that $\nu_2(\lambda_j+1)=a$ for all $\lambda_j\neq -1$.
\end{enumerate}
\end{lem}

\begin{proof}
We prove (1) by way of contradiction. Suppose there is an integer $a$ such that $\nu_2(\lambda_j)=a$ for all $1 \leq j \leq n.$ Let $\lambda'_j = \lambda_j/2^a.$ Then each $\lambda'_j$ is an odd integer. Since $n$ is odd, we get
    \[ \sum_{j=1}^n \lambda'_j  \equiv \sum_{j=1}^n  1 \equiv n  \equiv 1 \pmod{2}.\]
But since the trace of $A(\Gamma)$ is zero, we get
    \[ \sum_{j=1}^n \lambda'_j = \frac{1}{2^a} \sum_{j=1}^n \lambda_j= 0,\]
a contradiction. To prove (2), suppose that $\lambda_i$ is odd and there exists an integer $a$ such that $\nu_2(\lambda_j+1)=a$ for all $1 \leq j \leq n.$ Since $\lambda_i$ is odd, we know that $a \geq 1.$ Consequently $\lambda_j$ is odd for all $1 \leq j \leq n.$ Since $n$ is odd, we get that $\sum_{j=1}^n \lambda_j$ is odd, which is impossible because the trace of $A(\Gamma)$ is zero. Lastly, suppose $0$ is an eigenvalue with even multiplicity. If there is an $a\in\Z$ such that $\nu_2(\lambda_i)=a$ for all $\lambda_i\neq 0$, then a similar argument used in the proof of (1) yields the desired conclusion. Meanwhile, if $-1$ is an eigenvalue of $\Gamma$ and $\lambda_j\neq-1$ is odd for some $j$, then the same argument used in the proof of (2) yields the desired conclusion. 
\end{proof}

\begin{rem}
Lemma \ref{lem:congruence_condition_join}(2) is false if the condition that at least one eigenvalue is odd is not met. For example, the star graphs $K_{1, 4k^2}$ have an odd number of vertices and $\sigma(K_{1, 4k^2})=\{2k, -2k, 0\}$. So no eigenvalue of $K_{1, 4k^2}$ is odd and $\nu_2(\mu_i+1)=0$ for all $i$.
\end{rem}

\begin{cor}
Let $\Gamma$ be an integral graph where $|V(\Gamma)|$ is odd. If vertex $u$ has a false twin in $\Gamma$ such that $\Phi_u=\sigma(\Gamma)$ and $0\in\sigma(\Gamma)$ has even multiplicity, then $u$ is sedentary in $\Gamma$. If vertex $u$ has a true twin in $\Gamma$ such that $\Phi_u=\sigma(\Gamma)$ and some $\lambda\in\sigma(\Gamma)\backslash\{ -1 \}$ is odd, then $u$ is sedentary in $\Gamma$.
\end{cor}

\begin{proof}

If $u$ belongs to a false twin set of size at least three, then $u$ is sedentary by Theorem \ref{thm:twinschar}. Thus, we may assume that the twin set containing $u$ has size two. By assumption, $\Phi_u=\sigma(\Gamma)$ and $\theta=0$ has even multiplicity.
If $u$ admits PST in $\Gamma$, then Theorem \ref{thm:parcon} yields $\nu_2(\lambda-\theta)=a$ for all $\lambda\in\sigma(\Gamma)\backslash\{\theta\}$, a contradiction to Lemma \ref{lem:congruence_condition_join}(3). The same argument applies if $u$ has a true twin in $\Gamma$ with $\Phi_u=\sigma(\Gamma)$, $\theta=-1$,  and some $\lambda\in\sigma(\Gamma)\backslash\{\theta\}$ is odd. In both cases, $u$ is not involved in PST, and so it is sedentary in $\Gamma$ by Theorem \ref{thm:twdich}.
\end{proof}



We are now ready to prove the following result.


\begin{thm}
\label{thm:twlex}
Suppose $\Gamma=\Gamma_1[\Gamma_2]$ is an integral graph where $\Gamma_2\in\{O_2,K_2\}$. If $|V(\Gamma)|\equiv 2$ (mod 4) and $(u,j)\in V(\Gamma)$ such that $\Phi_{\be_{u}}(\Gamma_1)=\spec(\Gamma_1)$, then the following hold.
\begin{enumerate}
\item If $\Gamma=\Gamma_1[O_2]$, then vertex $(u,j)$ is sedentary in $\Gamma$.
\item If $\Gamma=\Gamma_1[K_2]$, then vertex $u$ is sedentary in $\Gamma$ if and only if either
(i) $-1\in\spec(\Gamma_1)$ or (ii) $-1\notin\spec(\Gamma_1)$ and some eigenvalue in $\spec(\Gamma_1)$ is odd.
\end{enumerate}
\end{thm}

\begin{proof}
Let $\Gamma=\Gamma_1[\Gamma_2]$, where $\Gamma_2\in\{O_2,K_2\}$. Since $|V(\Gamma)|\equiv 2$ (mod 4), it follows that $|V(\Gamma_1)|$ is odd. From (\ref{eq:lexprod}), $\spec(\Gamma_1[O_2])=\{2\mu:\mu\in\spec(\Gamma_1)\}\cup\{\theta=0\}$ and $\spec(\Gamma_1[K_2])=\{2\mu+1:\mu\in\spec(\Gamma_1)\}\cup\{\theta=-1\}$. Since $\sigma(\Gamma)\subseteq\Z$, it follows from the spectra of $\Gamma_1[O_2]$ and $\Gamma_1[K_2]$ that $\spec(\Gamma_1)\subseteq\Z$. Now, if (a) $\Gamma=\Gamma_1[O_2]$ and $0\in \spec(\Gamma_1)$ or (b) $\Gamma=\Gamma_1[K_2]$ and $-1\in \spec(\Gamma_1)$, then vertex $(u,j)$ is sedentary in $\Gamma$ by Corollary \ref{cor:blup}. So to prove (1), it suffices to consider the case when $0\notin \spec(\Gamma_1)$, and to prove (2), it suffices to consider the case when $-1\notin \spec(\Gamma_1)$ and some eigenvalue of $\spec(\Gamma_1)$ is odd. 

Assume that condition (2i) or (2ii) above holds. By way of contradiction, suppose $u$ is not sedentary in $\Gamma$. Since $\Phi_{\be_{(u,j)}}(\Gamma)\subseteq\Z$, this assumption is equivalent to vertex $u$ being involved in PST in $\Gamma$ by Theorem \ref{thm:twdich}. By Theorem~\ref{thm:parcon}, $\nu_2(\lambda-\theta)$ is a constant for all $\lambda\in\sigma(\Gamma)\backslash\{\theta\}$. Let us first deal with the case when $\Gamma=\Gamma_1[O_2]$, and let us assume $0\notin \sigma(\Gamma_1)$. The condition that $\nu_2(\lambda)$ is a constant for all $\lambda\in\sigma(\Gamma)\backslash\{0\}$ is equivalent to $\nu_2(\mu)$ is a constant for all $\mu\in\spec(\Gamma_1)$. Since $|V(\Gamma_1)|$ is odd, Lemma \ref{lem:congruence_condition_join}(1) yields a contradiction. Next, suppose $\Gamma=\Gamma_1[K_2]$ and assume $-1\notin \spec(\Gamma_1)$. The condition that $\nu_2(\lambda+1)$ is a constant for all $\lambda\in\sigma(\Gamma)\backslash\{-1\}$ is equivalent to $\nu_2(\mu+1)$ is a constant for all $\mu\in\sigma(\Gamma_1)$. If some eigenvalue in $\sigma(\Gamma_1)$ is odd, then because $|V(\Gamma_1)|$ is odd and $-1\notin \sigma(\Gamma_1)$, applying Lemma \ref{lem:congruence_condition_join}(2) again yields a contradiction. In both cases, vertex $u$ is not involved in PST, and hence it is sedentary by Theorem \ref{thm:twdich}. To complete the proof of (2), we now show that $\Gamma_1[K_2]$ is not sedentary whenever $-1\notin\spec(\Gamma_1)$ and $\spec(\Gamma_1)$ consists of all even integers. Indeed, if $\Gamma=\Gamma_1[K_2]$, then $\sigma(\Gamma)=\{2\mu+1:\mu\in\spec(\Gamma_1)\}\cup\{\theta=-1\}$, and as each $\mu\in\spec(\Gamma_1)$ is even, we get $\nu_2(2\mu+1-(-1))=1+\nu_2(\mu+1)=1$. Invoking Theorem \ref{thm:sedwreath} completes the proof.
\end{proof}

Recall that a subset $S$ of $G$ is called \textit{normal} if $S$ is stable under conjugation. In other words, $S$ is a union of some conjugacy classes in $G.$ We say that $\Gamma=\Gamma(G,S)$ is a normal Cayley graph if $S$ is a normal subset of $G$. 

\begin{expl}
Let $\Gamma=\Gamma(G,S)$ be an integral normal Cayley graph of odd order. Then $\Gamma$ has an odd eigenvalue \cite{arnadottir2025note}. Combining this with Theorems \ref{thm:twlex} and \ref{thm:twinschar}(1), we get that $\Gamma[O_m]$ and $\Gamma[K_m]$ are sedentary for all $m\geq 2$.
\end{expl}

\subsection{Cayley graphs}

Here, we investigate whether the condition that all eigenvalues of a Cayley graph $\Gamma$ are divisible by some prime $p$ has structural implications on $\Gamma$.

\begin{cor} \label{cor:2mod4}
Let $G$ be a group such that $|G| \equiv 2 \pmod{4}$, $|G|>2$,  and $G$ has a normal $2$-Sylow subgroup; say $P$. If $\Gamma(G,S)$ is an integral normal Cayley graph with twins such that $S  \not \subset P$, then $\Gamma(G,S)$ is sedentary. 
\end{cor}

\begin{proof}
Let  $|G| \equiv 2 \pmod{4}$ and $\Gamma=\Gamma(G,S)$ be an integral normal Cayley graph with twins. 
Suppose to the contrary that $\Gamma(G,S)$ is not sedentary. By  \Cref{prop:existence_of_twins}, we know that there exists a non-trivial subgroup $H$ of $G$ such that $H$ is a set of twins. Since $\Gamma(G,S)$ is not sedentary, we have $|H|=2$ by Theorem \ref{thm:caytwsed}. By assumption, 
 $|G|\equiv 2\pmod 4$ and $|H|=2$, we know that $H$ is a $2$-Sylow subgroup of $G.$ Since all $2$-Sylow subgroups are conjugate and $P$ is normal, we must have $P=H.$ By  \Cref{lem:conditions_for_twins} and \Cref{lem:conditions_for_twins2}, Remark \ref{rem:Gamma1} and (\ref{eq:lexprod}), we have that $\Gamma=\Gamma_1[\Gamma_2]$ where $\Gamma_2\in\{O_2,K_2\}$, $|\Gamma_1|$ is odd, and $\Gamma_1$ is a Cayley graph over $G/H$. More specifically, we have $\Gamma_1 \cong \Gamma(G/H, S')$ where $S'$ is the image of $S$ under the projection map $G \to G/H$ (excluding the coset $eH$.) Moreover, since $S$ is a union of conjugacy classes, so is $S'$. Additionally, since $S  \not \subset P$, $S'$ is not empty. This shows that $\Gamma_1$ is an integral normal Cayley graph over $G/H$. Since $|G/H|$ is odd and $S'$ is not empty,  $\Gamma_1$ has an odd eigenvalue by \cite{arnadottir2025note}. By \cref{thm:twlex}, we conclude that $\Gamma$ must be sedentary. 
\end{proof}

Let $\Gamma=\Gamma(G,S)$ be a integral Cayley graph and $H$ be a subgroup of $G$ that is a twin set in $\Gamma$ with $m=|H|$. By (\ref{eq:lexprod}), the eigenvalues of $\Gamma$ are multiples of $m$ if $\Gamma\cong \Gamma_1[O_m]$, and congruent to $-1$ modulo $m$ otherwise. This prompts us to ask: to what extent does the converse of this fact hold? That is, if all eigenvalues of $\Gamma$ are multiples of $m$, then does $\Gamma$ contain a false twin set of size $m$? Similarly, if all eigenvalues of $\Gamma$ are congruent to $-1$ modulo $m$, does $\Gamma$ contain a true twin set of size $m$? We answer these questions to the negative 
by characterizing the occurrence of this modular property on $\sigma(\Gamma)$ under mild assumptions on $S.$

\begin{thm} 
\label{thm:even-eig}
Let $G$ be a finite group, $p$ be a prime that divides $|G|$, and $P$ be a normal $p$-Sylow subgroup of $G$. Let $S\subset G$ be such that $\Gamma(G,S)$ is integral. Let
$A_S=\sum_{s\in S}s\in \Z[G]$ and let $\bar{A}_S$ be its image in $\F_p[G]$. The following statements hold.
\begin{enumerate}

\item If $\bar{A}_S$ belongs to the kernel of the natural map
$\Phi:\F_p[G]\to \F_p[G/P]$, then all eigenvalues of $\Gamma(G,S)$ are divisible by $p$.
\item If, in addition, $\Gamma(G,S)$ is normal, then the converse also holds. Equivalently, all eigenvalues of $\Gamma(G,S)$ are divisible by $p$ if and only if $\bar{A}_S\in\ker(\Phi)$.
\end{enumerate}
Moreover, the condition $\bar{A}_S\in\ker(\Phi)$ is equivalent to the condition that for every coset $gP$,
\[ |S\cap gP|\equiv 0\pmod p.\]
\end{thm}

\begin{proof}
Let $A$ be the adjacency matrix of $\Gamma(G,S)$. This is the matrix of the linear map
\[
\Q[G]\to \Q[G],\qquad g\mapsto gA_S.
\]
Let $f(x)=\prod_{i=1}^{|G|}(x-\lambda_i)$ be the characteristic polynomial of $\Gamma$. Since $P$ is a normal $p$-Sylow subgroup, the Jacobson radical of $\F_p[G]$ is the kernel of the following map 
\[
J(\F_p[G])=\ker\left(\Phi:\F_p[G]\to\F_p[G/P]\right)
\]
(see \cite[Theorem~16.6]{passman1971infinite}). If $\bar{A}_S\in\ker(\Phi)$, then $\bar{A}_S\in J(\F_p[G])$, and so $\bar{A}_S$ is nilpotent. Thus
\[
f(x)\equiv x^{|G|}\pmod{p},
\]
which implies that all the integer eigenvalues $\lambda_j$ are divisible by $p$.

Conversely, suppose that $\Gamma(G,S)$ is normal and all $\lambda_j$ are divisible by $p$. Then
\[
f(x)\equiv x^{|G|}\pmod{p}.
\]
By the Cayley--Hamilton theorem, $\bar{A}_S$ is nilpotent in $\F_p[G]$. Since $\Gamma(G,S)$ is normal, $S$ is a union of conjugacy classes. Thus, $S$ is stable under conjugation, and so $\bar{A}_S$ is central, implying that $\bar{A}_S\in J(\F_p[G])$. Thus $\bar{A}_S\in\ker(\Phi)$. Finally,
\[
\Phi(\bar{A}_S)
=
\sum_{gP\in G/P}|S\cap gP|\,gP,
\]
where the coefficients are taken modulo $p$. Hence $\Phi(\bar{A}_S)=0$ if and only if for every coset $gP$,
\[
|S\cap gP|\equiv 0\pmod{p}.\qedhere
\]
\end{proof}

We now derive corollaries of Theorem \ref{thm:even-eig}. 

\begin{cor}
\label{cor:pmod}
Let $p$ be a prime and $G$ be a finite group such that $\gcd(p, |G|)=1$. 
If $\Gamma(G,S)$ is an integral normal Cayley graph such that $S \neq \varnothing$, then $\Gamma(G,S)$ has at least one eigenvalue which is not a multiple of $p$.  
\end{cor}

\begin{proof}
In this case,  $P = \{1 \}$ is the only $p$-Sylow subgroup of $G.$ Since $S$ is nonempty, we can choose $s \in S.$ Then $|S \cap sP|=1$ which is not a multiple of $p.$ Therefore, by \Cref{thm:even-eig}, we know that $\Gamma(G,S)$ has an eigenvalue which is not a multiple of $p.$
\end{proof}

Taking $p=2$ in Corollary \ref{cor:pmod} recovers a result in \cite{arnadottir2025note}, which states that an integral normal Cayley graph on an odd number of vertices has an odd eigenvalue.

Next, we deal with the case where $p \mid |G|$ but $p^2 \nmid |G|.$ In this case, we give a partial answer to the question raised immediately after \Cref{cor:2mod4}. Namely, if all eigenvalues of $\Gamma(G,S)$ are multiples of $p$, does it imply that $\Gamma(G,S)$ has a set of twins of size $p$?
\begin{cor}
Let $G$ be a finite group and $p$ a prime number with $p \mid |G|$ and $p^2 \nmid |G|$. Suppose that $G$ has a normal $p$-Sylow subgroup $P$ such that $|P|=p$. If $\Gamma(G,S)$ is an integral normal Cayley graph, then the following hold. 
    \begin{enumerate}
    \item All eigenvalues of $\Gamma(G,S)$ are multiples of $p$ if and only if $\Gamma(G,S)$ is isomorphic to $\Gamma_1[O_p]$ where $\Gamma_1$ is a Cayley graph on $G/P.$
    \item All eigenvalues of $\Gamma(G,S)$ are congruent to $-1$ modulo $p$ if and only if $\Gamma(G,S)$ is isomorphic to $\Gamma_1[K_p]$ where $\Gamma_1$ is a Cayley graph on $G/P.$
    \end{enumerate}
\end{cor}

\begin{proof}

We first prove (1). First, if $\Gamma(G,S) = \Gamma_1[O_p]$, then from (\ref{eq:lexprod}) we know that
$$
\sigma\bigl(\Gamma_1[O_p]\bigr) = \{p\lambda:\lambda\in\sigma(\Gamma_1)\} \cup \{0\}\},
$$
This shows that all eigenvalues of $\Gamma$ are multiples of $p.$

Suppose now that all of the eigenvalues of $\Gamma(G,S)$ are multiples of $p$. Then by \Cref{thm:even-eig}, we know that $|S \cap gP|$ is a multiple of $p$ for all $g.$ In particular, $P \cap S = \varnothing$ since $1 \not \in S$ and for each $g \not \in P$, either $S \cap gP = \varnothing$ or $gP \subset S.$ In other words, $S$ is a union of cosets of $P$. By \Cref{lem:conditions_for_twins2}, we know that $\Gamma(G,S)$ has a false twin set of size $p$ and $\Gamma(G,S) \cong \Gamma_1[O_p]$ where $\Gamma_1$ is a graph on the coset $G/P.$ Since $P$ is normal, $G/P$ is a group and $\Gamma_1$ is a Cayley graph on $G/P.$ 
The proof of (2) follows similarly by considering $A_{S^*}$ where $S^*$ is the set formed by adding  the identity to $S$, in which case, the matrix representation of the action of $A_{S^*}$ on $\mathbb{Q}[G]$ is $A+I_n$.
\end{proof}

\begin{expl} \label{expl:D4}
Let $G=D_4$ be the dihedral group generated by $r,s$,
\[D_4=\langle r,s\mid r^4=s^2=1,\ srs=r^{-1}\rangle. \]
Let $\Gamma=\Gamma(G,S)$ where $S=\{r,r^2,r^3,s\}$. Then $\spec(\Gamma)=\{4^{(1)},2^{(1)},0^{(3)},(-2)^{(3)}\}$. As $G$ is a $2$-group, it is its own $2$-Sylow group. Since $|S|$ is even, \Cref{thm:even-eig}(1) guarantees that all eigenvalues of $\Gamma(G,S)$ are even. 
However, one verifies that $\Gamma$ has no twin set of size $2$ or $4$. 
\end{expl}

Finally, we deal with the case where $G$ is a $p$-group. In this case, $G$ is its own $p$-Sylow subgroup and we can remove the normality condition on $\Gamma(G,S)$.

\begin{cor}
Let $p$ be a prime and $G$ be a finite $p$-group. If $\Gamma(G,S)$ is integral, then
every eigenvalue $\lambda$ of $\Gamma(G,S)$ satisfies
$\lambda\equiv |S|\pmod p.$
In particular, all eigenvalues of $\Gamma(G,S)$ are divisible by $p$
if and only if $|S|$ is divisible by $p$.
\end{cor}

\begin{proof}
We modify the proof for \Cref{thm:even-eig}.  Let $A$ be the adjacency matrix of $\Gamma(G,S)$ and let $A_S=\sum_{s\in S}s\in\Z[G].$
The image of $A_S-|S|1_G$ in $\F_p[G]$ has augmentation zero. Since
$G$ is a finite $p$-group, the augmentation ideal of $\F_p[G]$ is
nilpotent. It follows that $A-|S|I$ is nilpotent modulo $p$, and hence
$f_A(x)\equiv (x-|S|)^{|G|}\pmod p.$  Since $\Gamma(G,S)$ is integral, we have $f_A(x)=\prod_{i=1}^{|G|}(x-\lambda_i)$ 
where $\lambda_1,\ldots,\lambda_{|G|}$ are the eigenvalues of $\Gamma(G,S).$ Therefore, $\lambda_i\equiv |S| \pmod{p}$ for every $i$. The second statement follows immediately.
\end{proof}

\medskip

\section{Future work}
\label{sec:futwork}

In this work, we investigated sedentariness in walk-regular graphs and Cayley graphs. We imported tools from algebraic number theory to characterize sedentariness in integral walk-regular graphs in terms of their necklace polynomials. As an application, we showed that an integral Cayley graph of odd order is sedentary. Under mild assumptions, we also characterized integral normal Cayley graphs whose every eigenvalue is divisible by a prime $p$. To inspire future work on this topic, we pose the following problems.

In \cite{bavsic2013characterization}, PST was characterized in integral Cayley graphs over $\Z_n$ (also called integral circulant graphs or ICGs). By Theorem~\ref{thm:oddsed}, ICGs with an odd number of vertices are sedentary. It would be interesting to characterize sedentariness in ICGs with an even number of vertices.

For sedentary vertex-transitive graphs with twins, we were able to provide a bound on their sedentariness whenever condition (1) or (2) of Theorem~\ref{thm:caytwsed} holds. However, if condition (3) or (4) of Theorem~\ref{thm:caytwsed} holds, then we currently do not have a bound for sedentariness. In this case, \cite[Corollary 3.14]{Monterde2022} implies that the twins are strongly cospectral. In general, it would be interesting to find a bound for sedentary twins that admit strong cospectrality.

We provided sedentary families of direct products of complete graphs. Our computations reveal that the Cartesian product and the direct product of the same set of complete graphs are either both sedentary, or both not sedentary. Is it the case that the Cartesian product of complete graphs is sedentary if and only if their direct product is? 

Finally, \Cref{thm:even-eig} raises an interesting question: can we classify integral Cayley graphs whose eigenvalues have equal $p$-adic valuation, for some fixed prime number $p$? 
Additionally, is there a result analogous to Theorem~\ref{thm:even-eig} for powers of $p$?  

\bigskip

\noindent \textbf{AI Usage Declaration}
We acknowledge the use of ChatGPT 5.6 in developing the proof of Proposition \ref{prop:bound-coefficients} and checking our examples. AI tools were not used to draft the manuscript. All formal statements, arguments, and proofs in the manuscript were written and checked by the authors, who take full responsibility for the accuracy and integrity of the article.

\bigskip

\noindent \textbf{Acknowledgments.}
KM is supported by NSERC Discovery Grant RGPIN-2025-04101. HM is supported by the Pacific Institute for the Mathematical Sciences. TTN is supported by an AMS-Simons Travel Grant. Parts of this work were done while TTN was visiting faculty at Argonne National Lab. He thanks Argonne for their hospitality.

\bigskip

\bibliographystyle{amsplain}
\bibliography{references.bib}
\end{document}